\documentclass[11pt]{amsart}
\usepackage{latexsym,xy}
\xyoption{all}

\newtheorem{theorem}{Theorem}[section]
\newtheorem{lemma}[theorem]{Lemma}
\newtheorem{proposition}[theorem]{Proposition}

\theoremstyle{definition}

\newtheorem{example}[theorem]{Example}
\newtheorem{remark}[theorem]{Remark}

\renewcommand{\P}{{\mathbb{P}}}
\newcommand{\Z}{{\mathbb{Z}}}
\newcommand{\A}{{\mathbb{A}}}

\newcommand{\Fc}{{\mathcal{F}}}
\newcommand{\Oc}{{\mathcal{O}}}
\newcommand{\syz}{{\mathrm{Syz}}}

\numberwithin{equation}{section}

\begin{document}

\title{A case study in bigraded commutative algebra} 

\maketitle

\section*{Introduction}

The purpose of this chapter is to illustrate how bigraded commutative
algebra differs from the classical graded case and to indicate some of
the theoretical tools needed to understand free resolutions in the
bigraded case.  Bigraded commutative algebra is a special case of
multigraded commutative algebra, which is in turn an instance of toric
(or polytopal) algebra.

Generalizing classical concepts of algebraic geometry or commutative
algebra to the multigraded setting is a very active area of current
research; a good overview of the area is given in \cite{sturm}.  For
instance, bigraded Castelnuovo-Mumford regularity is studied in
\cite{acd}, \cite{r}, and \cite{hw}, and generalizations to the
multigraded case are investigated in \cite{ms}.  Another important
topic is the study of the toric Hilbert scheme, which parameterizes
the set of all ideals having the same multigraded Hilbert function as
some fixed toric ideal.  References for this are \cite{ps1},
\cite{ps2}, \cite{mt}, and \cite{sst}.  In a different direction,
there is also the study of zero-dimensional schemes in the multigraded
setting, which began in \cite{gmr}, and is studied further in
\cite{hvt}, \cite{vt}, and \cite{svt}.

In this chapter we consider the bigraded commutative algebra of the
polynomial ring $R = k[x,y,z,w]$, where $k$ is algebraically closed,
and $x,y$ have degree $(1,0)$ and $z,w$ have degree $(0,1)$.  We will
write $R_{m,n}$ for the graded piece of $R$ in degree $(m,n)$.  Our
main object of study is the ideal $I \subset R$ generated by three
bihomogeneous polynomials
\[
p_{0}, p_{1}, p_{2} \in R_{2,1}
\]
assuming that the $p_i$ have no common zeros on $\P^{1} \times \P^{1}$.

In the singly graded case, it is well-known that if three homogeneous
polynomials in $k[s,t,u]$ have no common zeros on $\P^2$, then their
Koszul complex is a free resolution of the ideal they generate.  We
will see that life is more complicated (and more interesting!)\ in the
bigraded case.  In particular, for the above polynomials $p_{0},
p_{1}, p_{2}$, we will prove that
\begin{itemize}
\item the Koszul complex of $p_{0}, p_{1}, p_{2}$ is never exact, and
\item there are two possible minimal free resolutions of the ideal $I$.
\end{itemize}

In each section of the chapter, we highlight the tools from
commutative algebra that are used.  The techniques include regular
sequences, depth, spectral sequences, vanishing theorems,
determinants, the Segre embedding, resultants, and mapping cones.  Our
main results are the determination of the Hilbert function and minimal
free resolution of the ideal $I$.  In the last section of the chapter,
we outline some connections to the implicitization problem in
geometric modeling.
                                                     
\section{Preliminary Analysis of the Koszul Complex}
\label{preliminary}

\noindent {\bfseries \emph{Tools}:} Regular sequences, depth, the
Buchsbaum-Eisenbud criterion

\subsection{The Koszul Complex} In this section, we will assume that
$p_{0}, p_{1}, p_{2} \in R$ have degree $(d_1,d_2)$, $d_1,d_2 > 0$,
and no common zeros on $\P^{1} \times \P^{1}$.  This gives the Koszul
complex
\begin{equation}
\label{koszul}
\begin{aligned}
0 \rightarrow R(-3d_1,-3d_2) &\xrightarrow{\left[
\begin{matrix}
\!p_2\!  \\
\!-p_1\! \\
\!p_0\! 
\end{matrix} \right]} R(-2d_1,-2d_2)^3
\xrightarrow{\left[\begin{matrix}
\!p_1 & \!p_2&\!0\! \\
\!-p_0 & \!0 &\!p_2\!\\
\!0 & \!-p_0 & \!-p_1\! 
\end{matrix} \right]}\\
{}& R(-d_1,-d_2)^3
\xrightarrow{\left[ \begin{matrix}
p_0 & p_1& p_2
\end{matrix}\right]} R \rightarrow 0. 
\end{aligned}
\end{equation}

\subsection{The Failure of Exactness} Our first observation is that
\eqref{koszul} is not exact.  This follows by analyzing the zero locus
of $p_{0}, p_{1}, p_{2}$ in $\A^4$.  Since $d_1,d_2 > 0$, we know that
$p_{0}, p_{1}, p_{2}$ vanish on
\[
\mathbf{V}(xz,xw,yz,yw) = (\A^2 \times \{0\}) \cup (\{0\} \times \A^2)
\subset \A^4.
\]
Furthermore, $p_{0}, p_{1}, p_{2}$ can't vanish simultaneously outside
this variety since they have no common zeros on $\P^{1} \times
\P^{1}$.  Thus $I = \langle p_{0}, p_{1}, p_{2}\rangle$ has a
two-dimensional zero locus in $\A^4$, i.e, $\mathrm{codim}(I) = 2$. 

Now assume that the Koszul complex is exact. Then:
\begin{itemize}
\item $p_{0}, p_{1}, p_{2}$ form a regular sequence
  \cite[Thm.\ 17.6]{e}, so that $\mathrm{depth}(I) = 3$.  
\item $R$ is Cohen-Macaulay \cite[Prop.\ 18.9]{e}, so that
  $\mathrm{codim}(I) = \mathrm{depth}(I)$.
\end{itemize}
It follows that $\mathrm{codim}(I) = 3$, which contradicts the
previous paragraph.  Thus \eqref{koszul} is not exact.

\subsection{A More Careful Analysis} There is more to say, for
\eqref{koszul} still has some exactness.  Since $I$ is the image of
$R(-d_1,-d_2)^3 \to R$, we can replace \eqref{koszul} by 
\begin{equation}
\label{koszul2}
0 \to R(-3d_1,-3d_2) \to R(-2d_1,-2d_2)^3 \to R(-d_1,-d_2)^3 \to I \to
0.
\end{equation}

\begin{lemma}
\label{exex}
The complex \eqref{koszul2} is exact except at
$R(-d_1,-d_2)^3$. 
\end{lemma}

\begin{proof} 
We will give two proofs here and a third in the next section.  For the
first proof, we use a more detailed version of the relation between
depth and the Koszul complex.  Recall from above that
$\mathrm{depth}(I) = \mathrm{codim}(I) = 2$.  If we denote the Koszul
complex \eqref{koszul} by
\[
K: 0 \longrightarrow K_0 \xrightarrow{\varphi_0} K_1 \xrightarrow{\varphi_1}
K_2 \xrightarrow{\varphi_2} K_3 \longrightarrow 0,
\]
then \cite[Thm.\ 17.4]{e} tells us that
\[
H_j(K) = 0, j < r, \, \mathrm{and} \, H_r(K) \ne 0 \Longrightarrow \mathrm{depth}(I) =
r.
\]
In particular, we must have $H_0(K) = H_1(K) = 0$, since otherwise  $I$
would have depth $0$ or $1$.  Then we are done since
$R(-d_1,-d_2)^3 \to I$ is onto.

For second proof, we apply the Buchsbaum-Eisenbud criterion
\cite[Thm.\ 20.9]{e} to the complex 
\begin{equation}
\label{koszul3}
0 \longrightarrow K_0 \xrightarrow{\varphi_0} K_1 \xrightarrow{\varphi_1}
K_2.
\end{equation}
The Fitting ideals $I_i$ \cite[Sect.\ 20.2]{e} of $\varphi_0$ and
$\varphi_1$ are computed using the $i\times i$ minors of these maps.
One easily computes that
\[
I_1(\varphi_0) = I,\quad I_2(\varphi_0) = \{0\}
\]
and 
\[
I_2(\varphi_1) = I^2,\quad I_3(\varphi_1) = \{0\},
\]
so that $\varphi_0$ has rank $1$ and $\varphi_1$ has rank $2$.
Furthermore, $I^2$ and $I$ have the same depth since $\sqrt{I^2} =
\sqrt{I}$ \cite[Cor.\ 17.8]{e}.  Using $\mathrm{depth}(I) = 2$, we
obtain
\[
\mathrm{depth}(I_1(\varphi_0)) = \mathrm{depth}(I_2(\varphi_1)) = 2.
\]
Thus the hypotheses of the Buchsbaum-Eisenbud criterion are satisfied,
so that \eqref{koszul3} is exact.  As above, the lemma follows.
\end{proof}  

\section{Describing the Failure of Exactness}
\label{describing}

\noindent {\bfseries \emph{Tools}:} Sheaf cohomology, spectral
sequences, \v{C}ech cohomology

\subsection{The Syzygy Module} From now on, we will assume that
$p_{0}, p_{1}, p_{2} \in R$ have degree $(2,1)$ and don't vanish
simultaneously on $\P^{1} \times \P^{1}$.  Recall that the
\emph{syzygy module} $\syz(p) = \syz(p_0,p_1,p_2)$ is defined by the
exact sequence
\[
0 \longrightarrow \syz(p) \longrightarrow R(-2,-1)^3
\xrightarrow{\left[\begin{matrix} p_0 & p_1 & p_2
\end{matrix}\right]} I \longrightarrow 0 
\]
Then the Koszul complex gives a map
\[
R(-4,-2)^3 \longrightarrow \syz(p)
\]
which is not surjective by Lemma~\ref{exex}.  Syzygies in the image of
this map are called \emph{Koszul syzygies}, and those not in the image
are \emph{non-Koszul}.  

\subsection{The Sheaf-Theoretic Koszul Complex}
We will use sheaf cohomology and spectral sequences to give an
abstract description of the non-Koszul syzygies (explicit descriptions
will appear in later sections).  Our strategy will be to examine the
Koszul complex \eqref{koszul2} (with $(d_1,d_2) = (2,1)$) one degree
at a time.  In degree $(m,n)$, \eqref{koszul2} gives the complex
\begin{equation}
\label{kozr}
0 \longrightarrow R_{m-6,n-3} \longrightarrow R_{m-4,n-2}^3
\longrightarrow R_{m-2,n-1}^3 \longrightarrow I_{m,n} \longrightarrow 0,
\end{equation}
which is exact except at $R_{m-2,n-1}^3$ by Lemma~\ref{exex}. 

We begin with the sheaf-theoretic version of the Koszul complex.  Let
$\Oc := \Oc_{\P^{1} \times \P^{1}}$ and recall that $\mathcal{O}(m,n)$
is the line bundle on $\P^1\times\P^1$ with the property that
\[
H^0(\P^1\times\P^1,\Oc(m,n)) = R_{m,n}.
\]

The Koszul complex of $p_{0}, p_{1}, p_{2}$ is
\[
0 \to \mathcal{O}(-6,-3) \to \Oc(-4,-2)^{3} \to \mathcal{O}(-2,-1)^{3}
\to \Oc \to 0,
\]
so that tensoring by $\Oc(m,n)$ gives
\[
0 \! \to \mathcal{O}(m-6,n-3)\! \to \Oc(m-4,n-2)^{3}\! \to \Oc(m-2,n-1)^{3}\!
\to \Oc(m,n)\! \to 0.
\]
We write this complex as
\begin{equation}
\label{koszulsh}
0 \to \Fc_{0}(m,n) \to \Fc_{1}(m,n) \to \Fc_{2}(m,n) \to \Fc_{3}(m,n) \to 0.
\end{equation}
The key point is that \eqref{koszulsh} is exact.  To see why, recall
that for sheaves on $\P^{1} \times \P^{1}$, exactness is equivalent to
being exact at every point $P \in \P^{1} \times \P^{1}$.  Since
$p_{0}, p_{1}, p_{2}$ don't vanish simultaneously at $P$, they
generate the unit ideal in the local ring $\mathcal{O}_{\P^{1} \times
\P^{1},P}$.  By \cite[Prop.\ 17.14]{e}, it follows that the Koszul
complex is exact at $P$.

\subsection{The Spectral Sequences} As explained in \cite[pp.\
445--446]{gh}, the complex $\Fc_\bullet(m,n)$ has a hypercohomology
\[
\mathbb{H}^\bullet(\P^{1} \times \P^{1},\Fc_\bullet(m,n))
\]
and two spectral sequences
\begin{align*}
{}'\!E_2^{p,q}(m,n) &= H^p(\P^{1} \times
\P^{1},\mathcal{H}^q(\Fc_\bullet(m,n))) 
\Longrightarrow \mathbb{H}^\bullet(\P^{1} \times \P^{1},\Fc_\bullet(m,n))\\
{}''\!E_1^{p,q}(m,n) &= H^q(\P^{1} \times \P^{1},\Fc_p(m,n))
\qquad\,\Longrightarrow \mathbb{H}^\bullet(\P^{1} \times
\P^{1},\Fc_\bullet(m,n)).
\end{align*}
Since the Koszul complex \eqref{koszulsh} is exact, we have
$\mathcal{H}^q(\Fc_\bullet(m,n)) = 0$ for all $q$.  Thus the first
spectral sequence proves the standard fact that the hypercohomology of
a exact complex is trivial.  

It follows that the second spectral sequence becomes
\[
{}''\!E_1^{p,q}(m,n) = H^q(\P^{1} \times \P^{1},\Fc_p(m,n)) \Longrightarrow
0,
\]
where the differential $d_1$ is induced by the maps in the Koszul complex.

\subsection{Consequences of the Spectral Sequence} To examine this
spectral sequence, we set $H^{i}(m,n) := H^i(\P^{1} \times
\P^{1},\Oc(m,n))$.  Then the spectral sequence becomes
\[
\xymatrix@C=10.5pt{H^2(m-6,m-3) \ar[r]^{d_1} \ar[drr]^{d_2}
\ar[ddrrr]^{d_3} & H^2(m-4,n-2)^3\ar[drr]^{d_2} \ar[r]^{d_1} &
H^2(m-2,n-1)^3 \ar[r]^{\ \ \ \ \ d_1} & H^2(m,n) \\ H^1(m-6,m-3) \ar[r]^{d_1}
\ar[drr]^{d_2} & H^1(m-4,n-2)^3 \ar[r]^{\ d_1} \ar[drr]^{d_2} & H^1(m-2,n-1)^3
\ar[r]^{\ \ \ \ \ d_1} & H^1(m,n)\\ R_{m-6,n-3} \ar[r]^{d_1} & R_{m-4,n-2}^3
\ar[r]^{d_1} & R_{m-2,n-1}^3 \ar[r]^{d_1} & R_{m,n}}
\]
where the ${}''\!E_1^{p,q}(m,n)$ terms are shown.  This looks messy but
makes it easy to prove the following lemma.

\begin{lemma}
\label{syzE2}
For any $(m,n)$ we have an exact sequence
\[
0 \longrightarrow R_{m-6,n-3} \longrightarrow R_{m-4,n-2}^3
\longrightarrow Syz(p)_{m,n} \longrightarrow {}''\!E_2^{0,1}(m,n)
\longrightarrow 0
\]
where
\[
{}''\!E_2^{0,1}(m,n) = \mathrm{ker}\big(d_1\colon H^1(m-6,n-3)
\longrightarrow H^1(m-4,n-2)^3\big)
\]
and $d_1$ is induced by the corresponding map in the Koszul complex. 
\end{lemma}

\begin{proof}
Converging to $0$ means that ${}''\!E_\infty^{p,q}(m,n) = 0$ for all
$p,q$.  First consider $(p,q) = (1,0)$, which corresponds to
${}''\!E_1^{1,0}(m,n) = R_{m-4,n-2}^3$ in the above diagram.  All
differentials $d_r$, $r \ge 2$ that land at or originate from this
position must vanish.  It follows that
\[
0 = {}''\!E_\infty^{1,0}(m,n) = {}''\!E_2^{1,0}(m,n) =
\mathrm{ker}(d_1)/\mathrm{im}(d_1).  
\]
This proves exactness at $R_{m-4,n-2}^3$, and exactness at
$R_{m-6,n-3}$ is proved similarly.  (Note that this gives a third
proof of Lemma~\ref{exex}, as promised.)

Next consider ${}''\!E_1^{2,0}(m,n) = R_{m-2,n-1}$ and
${}''\!E_1^{0,1}(m,n) = H^1(m-6,n-3)$.  Here, we have to worry about
$d_2$, but all higher differentials vanish.  It follows that
\[
0 = {}''\!E_\infty^{2,0}(m,n) = {}''\!E_3^{2,0}(m,n)\quad \text{and}\quad  
0 = {}''\!E_\infty^{0,1}(m,n) = {}''\!E_3^{0,1}(m,n).
\]
However, the only way for these ${}''\!E_3$ terms to vanish is for the map
\begin{equation}
\label{d2iso}
d_2\colon {}''\!E_2^{0,1}(m,n) \longrightarrow {}''\!E_2^{2,0}(m,n)
\end{equation}
to be an isomorphism.  Since
\begin{align*}
{}''\!E_2^{0,1}(m,n) &= \mathrm{ker}(d_1)/\mathrm{im}(d_1) =
\mathrm{ker}(d_1)\\ 
{}''\!E_2^{2,0}(m,n) &= \mathrm{ker}(d_1)/\mathrm{im}(d_1) =
\syz(p)_{m,n}/\mathrm{im}(d_1), 
\end{align*}
the lemma follows immediately. 
\end{proof}

This lemma shows that ${}''\!E_2^{0,1}(m,n)$ gives a precise measure
of the non-Koszul syzygies.  In Section~\ref{measuring} we will use
vanishing theorems for sheaf cohomology to compute the size of
${}''\!E_2^{0,1}(m,n)$ for almost all $(m,n)$.

\subsection{Description of $d_2$} We will next unveil, in our special
case, the spectral sequence machinery that produces a non-Koszul
syzygy in degree $(m,n)$, out of an element $\varphi$ in the kernel of
the map $d_1\colon H^1(m-6,n-3) \to H^1(m-4, n-2)$.  The isomorphism
\[
d_{2}\colon \ker(d_1) \simeq \syz(p)_{m,n}/\mathrm{im}(d_1) \subset
R_{m-2,n-1}^{3}/\mathrm{im}(d_1)
\]
from \eqref{d2iso} can be explained as follows.

Consider a Leray covering $\mathcal U = (U_{i})$ of $\mathbb{P}^{1}
\times \mathbb{P}^{1},$ for instance the affine covering given by the
four open sets 
\begin{align*}
{}&U_1 =\{x \not=0, z \not =0\},\ U_2 =\{y \not=0, z
\not =0\}\\ {}&U_3 =\{x \not=0, w \not =0\},\ U_4 =\{y \not=0, w \not
=0\}.
\end{align*}
Then $H^1(m-6,n-3)$ is the \v{C}ech cohomology of this cover, so that
$\varphi$ is given by a collection of sections $\varphi_{ij}$ of
$\Oc(m-6,n-3)$ in $U_i \cap U_j$, with $\delta(\varphi_{ij})
=0.$ Since $d_1$ comes from the Koszul map, the condition
$d_1(\varphi) =0$ is equivalent to the existence of  sections
$(\alpha_i^0,\alpha_i^1,\alpha_i^2)$ of $\Oc(m-4,n-2)^3$ in each $U_i$, 
such that
\begin{equation}
\label{pkaik}
p_k \cdot \varphi_{ij} = \alpha_i^k - \alpha_j^k,
\quad \forall \, i,j,k.
\end{equation}
Then, since obviously $p_k \, p_\ell \, \varphi_{ij} =  p_\ell \,
p_k \, \varphi_{ij}$ for all $i,j,k, \ell$, one easily sees that
\[
p_\ell \, \alpha^k_i -  p_k \, \alpha^\ell_i =  p_\ell \,
\alpha^k_j - p_k \, \alpha^\ell_j, \quad \forall \, i,j,k.
\]
This defines a global section $\psi_{k\ell}$ of $\Oc(m-2,n-1)$.  Let
$\psi := (\psi_{12}, - \psi_{02}, \psi_{01})$. It is straightforward to
verify that this is a syzygy, i.e., that  the equality $p_{0} \psi_{12} - p_{1}
\psi_{02} + p_{2} \psi_{01} = 0$ holds.

However, the syzygy $\psi \in \syz(p)_{m,n}$ is not well-defined,
since it depends on the choice of $\alpha_i^k$ in \eqref{pkaik}.
Given another choice
\[
p_k \cdot \varphi_{ij}  =  \beta_i^k - \beta_j^k,
\quad \forall \, i,j,k,
\]
we get $\psi' \in \syz(p)_{m,n}$.  Comparing this equation to
\eqref{pkaik}, we see that $\beta_i^k - \alpha_i^k = \beta_j^k -
\alpha_j^k$ for all $i,j,k$.  Thus we have a global section $\gamma^k
= \beta_i^k - \alpha_i^k$ of $\Oc(m-4,n-2)$.  We leave it as an
exercise for the reader to verify that $\psi' = \psi + d_1(\gamma)$,
where $\gamma := (-\gamma^2,\gamma^1,-\gamma^0)$.  It follows that the
class $[\psi]$ of $\psi$ modulo the image of $d_1$ is well-defined.  A
more substantial exercise is to prove that $d_2(\varphi) = [\psi]$
using the \v{C}ech description of the spectral sequence of a complex
given in \cite[pp.\ 445--446]{gh}.

\section{A Picture of the Syzygy Module}
\label{measuring}

\noindent {\bfseries \emph{Tools}:} the K\"unneth formula, the
cohomology of line bundles on $\P^1$

\subsection{The Degree of a Syzygy} We defined $\syz(p)$ to be the
kernel of the map 
\[
 R(-2,-1)^3
\xrightarrow{\left[\begin{matrix} p_0 & p_1 & p_2
\end{matrix}\right]} I.
\]
It follows that
\[
\syz(p)_{m,n} \subset R(-2,-1)^3_{m,n} = R^3_{m-2,n-1}.
\]
This shows that syzygies in $\syz(p)_{m,n}$ are represented by vectors
of polynomials of degree $(m-2,n-1)$.  It will be important to keep
this in mind in the discussion that follows.

\subsection{Vanishing Cohomology} In Lemma~\ref{syzE2}, we saw that
the non-Koszul syzygies are determined by 
\begin{equation}
\label{e2for}
{}''\!E_{2}^{0,1}(m,n) = \ker(d_{1}\colon H^{1}(m-6,n-3) \to
H^{1}(m-4,n-2)^{3}). 
\end{equation}
In order to compute this, we first need to compute $H^1(m,n)$.  We do
this as follows.

\begin{lemma}
\label{vanishing}
If either $m,n \ge -1$ or $m,n \le -1$, then $H^1(m,n) = 0$.
\end{lemma}

\begin{proof}
Let $\pi_i\colon \P^1 \times \P^1 \to \P^1$ be projection onto the $i$th
factor.  Then 
\[
\Oc(m,n) = \pi_1^*\Oc_{\P^1}(m) \times
\pi_2^*\Oc_{\P^1}(n).  
\]
The K\"unneth formula \cite[Thm.\ 6.7.8]{ega} implies that
$ H^1(m,n)$ is isomorphic to
\[ H^0(\P^1,\Oc_{\P^1}(m))\otimes
H^1(\P^1,\Oc_{\P^1}(n)) \oplus
H^1(\P^1,\Oc_{\P^1}(m))\otimes H^0(\P^1,\Oc_{\P^1}(n)).
\]
Furthermore, we know that $H^0(\P^1,\Oc_{\P^1}(k)) = 0$ for $k
\le -1$, and the perfect pairing
\begin{equation}
\label{sd}
H^0(\P^1,\Oc_{\P^1}(-k-2))\times
H^1(\P^1,\Oc_{\P^1}(k)) \longrightarrow
H^1(\P^1,\Oc_{\P^1}(-2)) \simeq k
\end{equation}
described in \cite[Thm.\ 5.1 of Ch.\ III]{h} implies that
$H^1(\P^1,\Oc_{\P^1}(k)) = 0$ for $k \ge -1$.  The lemma now
follows easily.
\end{proof}

Note that \eqref{sd} is Serre Duality for $\P^1$ (see \cite[Thm.\ 7.1
of Ch.\ III]{h}).  We next apply the lemma to
${}''\!E_{2}^{0,1}(m,n)$.

\begin{proposition}
\label{easypart}
If $m \ge 5, n \ge 2$ or $m \le 5, n \le 2$, then
${}''\!E_{2}^{0,1}(m,n) = 0$.  Furthermore,
\begin{equation}
\label{m34n1}
\begin{aligned}
\dim_k{}''\!E_{2}^{0,1}(m,1) = \dim_k H^1(m-6,-2) &= \begin{cases} m-5
  & \ \ \, m \ge 6 \\ 0 & \ \ \, m \le 5\end{cases}\\
\dim_k {}''\!E_{2}^{0,1}(3,n) = \dim_k H^{1}(-3,n-3) &= \begin{cases}
  2(n-2) & \!n \ge 3 \\ 0 & \!n \le 2\end{cases}\\
\dim_k {}''\!E_{2}^{0,1}(4,n) = \dim_k H^{1}(-2,n-3) &= \begin{cases}
  n-2   & \ \ \ n \ge 3 \\ 0 & \ \ \ n \le 2.\end{cases}
\end{aligned}
\end{equation}
\end{proposition}

\begin{proof}
If $m \ge 5, n \ge 2$ or $m \le 5, n \le 2$, then $H^1(m-6,n-3) = 0$
by Lemma~\ref{vanishing}.  Then ${}''\!E_{2}^{0,1}(m,n) = 0$ follows
by \eqref{e2for}.  To study ${}''\!E_{2}^{0,1}(m,1)$, first observe
that $H^1(m-4,-1) = 0$ for all $m$ by Lemma~\ref{vanishing}.  Hence
\[
\dim_k{}''\!E_{2}^{0,1}(m,1) = \dim_k H^1(m-6,-2) = \dim_k
H^0(\P^1,\Oc_{\P^1}(m-6)),
\]
where the first equality uses \eqref{e2for} and the second uses the
K\"unneth formula and $H^1(\P^1,\Oc_{\P^1}(-2)) \simeq k$.  Then we
are done since 
\[
\dim_k H^0(\P^1,\Oc_{\P^1}(m-6)) \simeq \dim_k
k[x,y]_{m-6} = m-5.
\]
The other assertions are similar though slightly more complicated.  We
leave the details as an exercise.
\end{proof}

\subsection{The Picture}

For a given graded piece $\syz(p)_{m,n}$ of the syzygy module,
there are three things that could happen:
\begin{itemize}
\item[\LARGE{--}] the graded piece could be zero, or
\item[\LARGE{--}] it could be nonzero but consist only of Koszul
syzygies, or 
\item[\LARGE{--}] it could contain a non-Koszul syzygy.
\end{itemize}
\noindent Because $p_0, p_1,p_2$ have no common zeros on $\P^1 \times
\P^1$, we can predict in advance which of these possibilities occurs
for $\syz(p)_{m,n}$ for most $(m,n)$.  Here is the precise result.

\begin{proposition}
\label{syzpic}
The graded pieces $\syz(p)_{m,n}$ of the syzygy module can be
described as follows:
\begin{equation}
\label{mnpict}
\begin{matrix} 
n\  & \vdots & \vdots & \vdots & \vdots & \vdots & 
\vdots & \vdots & \vdots & \\[6pt]
5\ & 0 & 0 & ? & \Box & \Box & \bullet & \bullet & \bullet 
&\cdots\\[3pt]
4\ & 0 & 0 & ? & \Box & \Box & \bullet & \bullet & \bullet 
&\cdots\\[3pt]
3\ & 0 & 0 & ? & \Box & \Box & \bullet & \bullet & \bullet 
&\cdots\\[3pt]
2\ & 0 & 0 & 0 & 0 & \bullet & \bullet & \bullet & 
\bullet &\cdots\\[3pt]
1\ & 0 & 0 & 0 & 0 & 0 & 0 & \Box & 
\Box &\cdots\\[3pt]
0\ & 0 & 0 & 0 & 0 & 0 & 0 & 0 & 0 
&\cdots\\[6pt]
& 0 & 1 & 2 & 3 & 4 & 5 & 6 & 7 & m
\end{matrix}
\end{equation}
where in position $(m,n)$,
\begin{itemize}
\item[$0$] means that $\syz(p)_{m,n} = 0$.
\item[$\bullet$] means that $\syz(p)_{m,n} \ne 0$ and is generated
by the Koszul syzygies.
\item[$\Box$] means that $\syz(p)_{m,n}$ contains non-Koszul
syzygies.   The exact number is computed in Proposition~\ref{easypart}.
\item[$?$] means that we don't yet know what is going on.
\end{itemize}
\end{proposition}

\begin{proof}
When $m \le 1$ or $n = 0$, we clearly get zero since elements of
$\syz(p)_{m,n}$ are represented by polynomials of degree $(m-2,n-1)$.

If $m = 3$ or $n = 1$, then Lemma~\ref{syzE2} implies that
$\syz(p)_{m,n} = {}''\!E_2^{0,1}(m,n)$ for these degrees.  The
corresponding zeros and boxes follow from \eqref{m34n1}.

For $m = 4$, \eqref{m34n1} implies that ${}''\!E_2^{0,1}(4,2) = 0$ and
${}''\!E_2^{0,1}(4,n) \ne 0$ for $n \ge 3$.  Then Lemma~\ref{syzE2}
explains the bullet at position $(4,2)$ and the boxes that lie above
it.

Proposition~\ref{easypart} tells us that ${}''\!E_2^{0,1}(m,n) = 0$ when
$m \ge 5$ and $n \ge 2$.  Then the corresponding bullets follow from
Lemma~\ref{syzE2}.

Finally, ${}''\!E_2^{0,1}(2,2) = 0$ by Proposition~\ref{easypart}.
But the proposition gives no information about ${}''\!E_2^{0,1}(2,n) =
0$ for $n \ge 3$.  Hence the question marks.
\end{proof}

In \eqref{mnpict}, the generators of the Koszul syzygies lie in the
bullet closest to the origin, at position $(4,2)$.  To get some
insight into the question marks, consider the lowest question mark, at
position $(2,3)$.  Then $\syz(p)_{2,3} = {}''\!E_2^{0,1}(2,3)$ since
there are no Koszul syzygies in this degree.  One easily computes that
\[
\dim_k H^1(-4,0) = 3\quad\text{and}\quad \dim_k H^1(-2,1) = 2.
\]
Thus $\syz(p)_{2,3} = {}''\!E_2^{0,1}(2,3) = \mathrm{ker}(d_2
\colon H^1(-4,0) \to H^1(-2,1)^3)$ is the kernel of a map from $k^3$ to
$k^6$.  As we will soon see, there are two possibilities for its
dimension.

\subsection{Exactness Modulo $B$-Torsion} Let $B := \langle
xz,xw,yz,yw\rangle$.  If we think of the Koszul complex of sheaves on
$\A^4$, then we get
\[
0 \longrightarrow \Oc_{\A^4} \longrightarrow \Oc^3_{\A^4}
\longrightarrow \Oc^3_{\A^4} \longrightarrow \Oc_{\A^4}
\longrightarrow 0,
\]
where the maps are given by the matrices in \eqref{koszul}.  We showed
in Section~\ref{preliminary} that $\mathbf{V}(p_0,p_1,p_2) =
\mathbf{V}(B)$ since $p_0,p_1,p_2$ have no common zeros on
$\P^1\times\P^1$.  Hence the above complex is exact 
outside $\mathbf{V}(B)$.  

It follows that when we take global sections, the resulting sequence
is exact modulo $B$-torsion.  (We leave the proof as an exercise for
the reader.)  Since the sequence of global sections is the Koszul
complex \eqref{koszul}, we conclude that the Koszul complex is exact
modulo $B$-torsion.

In particular, every syzygy becomes Koszul by multiplying by elements
of $B^\ell$ for sufficiently large $\ell$.  Since the generators of
$B$ have degree $(1,1)$, this is easy to see in the picture
\eqref{mnpict}, for multiplying by an element of degree $(1,1)$ moves
up and to the right one unit.  Looking at the picture, it should be
clear that multiplying by $B^3$ makes every syzygy Koszul.

Finally, we note that the idea of being exact up to $B$-torsion plays
an important role in the concept of toric regularity described in
\cite{ms}.

\section{The Question Marks}
\label{QM}

\noindent {\bfseries \emph{Tools}:} the Segre embedding, planes in
$\P^5$, line bundles on $\P^1$ 

\subsection{The Segre Embedding}  Before we can study the question
marks in \eqref{mnpict}, we need to set up some notation.  Elements of
$R_{2,1}$ will be written
\[
a\,x^2z + b\,xyz + c\,y^2z + d\,x^2w + e\,xyw + f\,y^2w.
\]
We will identify $\P(R_{2,1})$ with $\P^5$.  Since $p_0,p_1,p_2 \in
R_{2,1}$ are linearly independent (otherwise they would have a common
zero on $\P^1\times\P^1$), they determine a plane
\[
W(p) := \P(\mathrm{Span}(p_0,p_1,p_2)) \subset \P^5.
\]

Reducible polynomials in $R_{2,1}$ can factor in two ways according to
the multiplication maps
\[
R_{2,0} \times R_{0,1} \to R_{2,1}\quad\text{and}\quad R_{1,1} \times
R_{1,0} \to R_{2,1}.
\]
The first will prove to be relevant to our situation.  It gives the
Segre embedding
\[
\P^2 \times \P^1 = \P(R_{2,0}) \times \P(R_{0,1}) \to \P^5 =
\P(R_{2,1}).
\] 
We will denote the image of this map as 
\[
Y \subset \P^5.
\]
It is well-known (and easy to verify using Macaulay 2) that $Y \subset
\P^5$ has degree 3, dimension 3, and is defined by the vanishing of the
$2\times2$ minors of the matrix
\begin{equation}
\label{23matrix}
\left[\begin{matrix} a & b & c \\ d & e & f \end{matrix}\right],
\end{equation}
where $a,\dots,f$ are homogeneous coordinates of $\P^5$.  This is
clear if we write the elements of $R_{2,1}$ in the form
\[
(a\,x^2 + b\,xy + c\,y^2)z + (d\,x^2 + e\,xy + f\,y^2)w,
\]
since the vanishing of the maximal minors of the matrix (\ref{23matrix})
is equivalent to the quadratic forms $a\,x^2 + b\,xy + c\,y^2$ and
$d\,x^2 + e\,xy + f\,y^2$ being proportional.

\subsection{The Two Cases} Our assumption that $p_0,p_1,p_2 \in
R_{2,1}$ have no common zeros on $\P^1 \times \P^1$ is very strong.
However, as we will now prove, it is not strong enough to determine
the question marks in \eqref{mnpict}.

\begin{theorem}
\label{th:main}
Let $W(p)$ and $Y$ be as above.  Then:
\begin{enumerate}
\item If $W(p)\cap Y$ is finite, then $\syz(p)_{2,n} = 0$ for all $n \ge
3$. 
\item If $W(p)\cap Y$ is infinite, then $\syz(p)_{2,n}$ has dimension
$n-2$ for all $n \ge 3$.
\end{enumerate}
\end{theorem}

\begin{remark}
Since $Y \subset \P^5$ has codimension 2 and $W(p)$ is a plane, $W(p)
\cap Y$ is finite when $W(p)$ is generic.  This observation leads to
the following terminology:
\begin{itemize}
\item $p_0,p_1,p_2$ are \emph{generic} if $W(p)\cap Y$ is finite.  In
this case, the question marks in \eqref{mnpict} are all zeros.
\item $p_0,p_1,p_2$ are \emph{non-generic} if $W(p)\cap Y$ is
infinite.  In this case, the question marks are all boxes (since
syzygies in $\syz(p)_{2,n}$ are non-Koszul).
\end{itemize}
The proof given below will show that in the non-generic case, $W(p) \cap
Y$ is a smooth conic in $W$.  We will also see that a suitable
parametrization of this conic gives an element of $\syz(p)_{2,3}$
whose multiples generate $\syz(p)_{2,n}$ for all $n \ge 3$.
\end{remark}

\begin{proof}
The graded piece $\syz(p)_{2,n}$ consists of syzygies $A_0 p_0 + A_1
p_1 + A_2 p_2 = 0$, where $A_0,A_1,A_2$ have degree $(0,n-1)$.
Suppose that $\syz(p)_{2,n} \ne 0$ for some $n \ge 3$.  Pulling out
the GCD, we can write this syzygy as the triple
\begin{equation}
\label{haaa}
h(z,w)(A_0(z,w),A_1(z,w),A_2(z,w)),
\end{equation}
where $A_0,A_1,A_2$ are coprime.
It is easy to see that $(A_0,A_1,A_2)$ is also a syzygy on
$p_0,p_1,p_2$.  
Furthermore,  $A_0,A_1,A_{2}$ are nonconstant by \eqref{mnpict}.

Now consider the map $\varphi\colon \P^1 \to W(p)$ defined by
\[
(s,t) \mapsto p_{s,t} := A_0(s,t)p_0 + A_1(s,t)p_1 + A_2(s,t)p_2.
\]
Since $A_0,A_1,A_2$ are nonconstant with $\mathrm{GCD}(A_0,A_1,A_2) =
1$, this is defined on all of $\P^1$ and the image is a curve.

However, $p_{s,t}(x,y,z,w)$ has the interesting property that
$p_{s,t}(x,y,s,t)$ is given by
\[
A_0(s,t)p_0(x,y,s,t) + A_1(s,t)p_1(x,y,s,t) +
A_2(s,t)p_2(x,y,s,t),
\]
which vanishes identically since $(A_0,A_1,A_2)$ is a syzygy on
$p_0,p_1,p_2$.  Fixing $s,t$, we have $p_{s,t}(x,y,s,t) = 0$ for all
$x,y$.  This easily implies that for this $s,t$, we have a
factorization
\[
p_{s,t}(x,y,z,w) = q(x,y)(tz-sw)
\]
for some $q(x,y) \in R_{2,0}$ (which depends on our fixed $s,t$).  This
factorization shows that in $\P^5$, $p_{s,t}$ gives a point of the
Segre variety $Y$.  It follows that $\varphi$ gives a map
\[
\varphi\colon \P^1 \to W(p)\cap Y.
\]
This proves that $W(p)\cap Y$ is infinite.  Part (1) of the theorem
follows immediately.

It remains to consider part (2) of the theorem.  For this, we need to
study $W(p) \cap Y$ when the intersection is infinite.  If $W(p) \cap Y$ has
dimension $> 1$, then we must have $W(p) \subset Y$ since $W(p)$ is a
plane.  This gives
\[
\P^2 \simeq W(p) \subset Y \simeq \P^2 \times \P^1.
\]
Since all morphisms $\P^2 \to \P^1$ are constant, it follows that
$W(p)$ corresponds to $\P^2 \times \{\mathrm{(c_1,c_2)}\} \subset \P^2
\times \P^1 \simeq Y$.  When we think of $Y$ as consisting of
reducible polynomials in $R_{2,1}$, we see that every polynomial in
$W(p)$ is divisible by the fixed linear form $\ell = c_1 z + c_2 w$.
This is impossible since $p_0,p_1,p_2$ have no common zeros on $\P^1
\times \P^1$.  Hence $\dim (W(p) \cap Y) > 1$ can't occur.

For the rest of the proof we will assume that $\dim (W(p) \cap Y) =
1$.  We first show that in this case, $W(p) \cap Y$ is contained in a
conic in $W(p) \simeq \P^2$.  To see why, recall that $Y = Q_1 \cap
Q_2 \cap Q_3$, where $Q_1,Q_2,Q_3 \subset \P^5$ are the quadric
hypersurfaces defined by the $2\times 2$ minors of \eqref{23matrix}.
By the above paragraph, we know that $W(p) \not\subset Y$, so that
$W(p) \not\subset Q_i$ for at least one $i$.  For this $i$, we have
\[
W(p) \cap Y \subset W(p) \cap Q_i,
\]
where $W (p)\cap Q_i$ is a conic in $W(p)$.  Thus $W(p) \cap Y$ 
lies in a conic in $W(p)$.

To complete the proof, we will show the following:
\begin{itemize}
\item[(a)] $W(p) \cap Y$ is a smooth conic in $W(p)$.
\item[(b)] A suitable parametrization of this conic gives a nonzero
element of $\syz(p)_{2,3}$.
\item[(c)] For $n \ge 3$, every element of $\syz(p)_{2,n}$ is
obtained from the element of $\syz(p)_{2,3}$ from {(b)} by
multiplication by a unique polynomial in $R_{0,n-3}$.
\end{itemize}
\noindent Note that part (2) of the theorem follows immediately from
(b) and (c).

We will prove (a), (b), and (c) simultaneously by studying curves
in $W(p) \cap Y$ parametrized by $\P^1$.  Suppose that
\begin{equation}
\label{param}
\phi \colon \P^1 \to W(p) \cap Y,
\end{equation}
is a parametrization which, when regarded as a map to $\P^5$, is given
by polynomials of degree $N$.  This means that $\phi^*
\Oc_{\P^5}(1)$ is $\Oc_{\P^1}(N)$.  However, we can also
view $\phi$ as a map to $Y \simeq \P^2 \times \P^1$ given by
\begin{equation}
\label{higi}
\phi(s,t) = (r_0(s,t),r_1(s,t),r_2(s,t))\times (h_0(s,t),h_1(s,t))
\in \P^2 \times \P^1,
\end{equation}
where $r_0,r_1,r_2$ have degree $\alpha$ and are relatively prime, and
$h_0,h_1$ have degree $\beta$ and are relatively prime.  It follows
that 
\[
\phi^*\Oc_{\P^2\times\P^1}(1,1) \simeq
\Oc_{\P^1}(\alpha+\beta).  
\]
Since $\Oc_{\P^5}(1)$ pulls back to
$\Oc_{\P^2\times\P^1}(1,1)$ under the Segre embedding, it
follows that $\Oc_{\P^5}(1)$ pulls back to
$\Oc_{\P^1}(\alpha+\beta)$ under $\phi$.  Comparing this
to our earlier description of the pullback, we conclude that
\[
N = \alpha+\beta.
\]

We next show that neither $\alpha$ nor $\beta$ can be zero.  Suppose,
for example, that $\alpha = 0$.  This implies that in $\P(R_{2,1})$,
we have
\[
\phi(s,t) = g(x,y)(h_0(s,t)z + h_1(s,t)w).
\]
where $g(x,y) \in R_{2,0}$ does not depend on $s,t$.  This gives a
two-dimensional subspace of $R_{2,1}$ whose elements are all divisible
by $g(x,y)$.  However, $\phi$ parametrizes a curve in $W(p) =
\P(\widetilde{W}(p))$, where $\widetilde{W}(p) =
\mathrm{Span}(p_0,p_1,p_2) \subset R_{2,1}$.  It follows that
$\widetilde{W}(p)$ has a two-dimensional subspace whose elements are
all divisible by $g(x,y)$.  Since $p_0,p_1,p_2$ have no common zeros
on $\P^1 \times \P^1$, an easy argument shows that any two linearly
independent elements of $\widetilde{W}(p)$ must be relatively prime.
This shows that the two-dimensional subspace constructed above can't
exist.  This contradiction proves that $\alpha > 0$, and $\beta > 0$
follows similarly.

Now we can prove (a).  We showed above that $W(p) \cap Y$ is
contained in a conic in $W(p)$.  If the conic is not smooth, then it
is a union of two lines, so that $W(p) \cap Y$ contains a line $\ell$.
This gives a parametrization \eqref{param} with $N = 1$.  Since $N =
\alpha+\beta$ with $\alpha$ and $\beta$ positive, this is impossible.
It follows that the conic must be smooth and hence must equal $W(p)
\cap Y$.  Thus (a) is proved.

We now turn to (b).  Since $W(p) \cap Y$ is a smooth conic, we can
find a parametrization \eqref{param} with $N = 2$.  Since $N =
\alpha+\beta$ with $\alpha, \beta > 0$, we must have $\alpha = \beta =
1$.  In $\P(R_{2,1})$, \eqref{higi} can be written
\begin{equation} \label{eq:phi}
\phi(s,t) = (r_0(s,t)x^2 + r_1(s,t)xy + r_2(s,t)y^2)(h_0(s,t)z +
h_1(s,t)w),
\end{equation}
where $r_0,r_1,r_2,h_0,h_1$ have degree 1 in $s,t$.  We also know that
$h_0,h_1$ are relatively prime, which means that they are linearly
independent in $R_{0,1}$.  Thus, by reparametrizing $\phi$, we may
assume that $h_0(s,t) = t$ and $h_1(s,t) = -s$.  Hence
\begin{equation}
\label{vpfor}
\phi(s,t) = (r_0(s,t)x^2 + r_1(s,t)xy + r_2(s,t)y^2)(tz - sw)
\end{equation}
in $\P(R_{2,1})$.  Since $\phi$ gives a curve in $W$, we can also
write $\phi(s,t)$ in the form
\begin{equation}
\label{vpfor2}
B_0(s,t)p_0(x,y,z,w) + B_1(s,t)p_1(x,y,z,w)
+ B_2(s,t)p_2(x,y,z,w).
\end{equation}
These expressions for $\phi(s,t)$ are the same.  The first vanishes
when $(s,t) = (z,w)$, so that the same is true for the second.  Hence
\[
B_0(z,w)p_0(x,y,z,w) + B_1(z,w)p_1(x,y,z,w)
+ B_2(z,w)p_2(x,y,z,w) = 0.
\]
Since $B_0(z,w),B_1(z,w),B_2(z,w)$ have degree 2 in $z,w$, it follows
that we get a nonzero element of $\syz(p)_{2,3}$.  This proves (b).

For (c), suppose that we have a nonzero syzygy in $\syz(p)_{2,n}$.  As
in the proof of part (1) of the theorem, we can write the syzygy in
the form \eqref{haaa}, namely
\[
h(z,w)(A_0(z,w),A_1(z,w),A_2(z,w)),
\]
where $A_0,A_1,A_2$ are relatively prime.  This gives a map $\P^{1}
\to W(p)\cap Y$ defined by
\[
(s,t) \mapsto p_{s,t}(x,y,z,w) := A_{0}(s,t)p_{0} + A_{1}(s,t)p_{1} +
A_{2}(s,t) p_{2}
\]
with the property that
\begin{equation}
\label{exprop}
p_{s,t}(x,y,s,t) = 0.
\end{equation}
The image is the smooth conic $W(p)\cap Y$.  Since $A_0,A_1,A_2$ are
relatively prime and $W(p)\cap Y$ is parametrized by $B_0,B_1,B_2$
from \eqref{vpfor} and \eqref{vpfor2}, one can show that there are
relatively prime polynomials $u(s,t),v(s,t)$ such that
\[
A_i(s,t) = B_i(u(s,t),v(s,t)),\quad i = 0,1,2.
\]
(We leave the details to the reader.)  Combining this 
with \eqref{vpfor} shows that 
\begin{align*}
p_{s,t} =\ &(r_0(u(s,t),v(s,t))x^2 + r_1(u(s,t),v(s,t))xy +
r_2(u(s,t),v(s,t))y^2)\times\\ &(v(s,t)z - u(s,t)w).
\end{align*}
Then \eqref{exprop} tells us that
\[
v(s,t)s - u(s,t)t = 0.
\]
This equation easily implies that
\[
u(s,t) = s\,g(s,t)\quad\text{and}\quad v(s,t) = t\,g(s,t).
\]
Since $u(s,t)$ and $v(s,t)$ are relatively prime, it follows that
$g(s,t)$ is a nonzero constant $c$.  This in turn implies that
\[
p_{s,t} = (r_0(cs,ct)x^2 + r_1(cs,ct)xy +
r_2(cs,ct)y^2)(ctz - csw) = c^{2}\phi(s,t).
\]
Hence $(A_0,A_1,A_2)$ is a constant times $(B_0,B_1,B_2)$.  Then (c)
follows.
\end{proof}

One comment is that if $p_0,p_1,p_2 \in R_{2,1}$ are allowed to be
arbitrary, then $W(p)$ can be an arbitrary plane in $\P^5$.  In
particular, $W(p) \cap Y$ can be a singular conic or a double line.
But once we require that $p_0,p_1,p_2$ have no common zeros on $\P^1
\times \P^1$, then $W(p) \cap Y$ is restricted to being finite or a
smooth conic.

\subsection{An Example} Here is an example to show that the
non-generic behavior described in Theorem~\ref{th:main} can actually
occur.

\begin{example} 
\label{ex:special}
Suppose that $p_{0} = x^{2}z$, $p_{1} = y^{2}w$ and $p_{2} = x^{2}w + 
y^{2}z$.  One easily checks that these polynomials have no common 
zeros on $\P^{1} \times \P^{1}$.  It is also easy to verify that
\[
w^{2}\, p_{0} + z^{2} \, p_{1} - zw \, p_{2} = 0.
\]
Thus $(w^2,z^2,-zw)$ gives a nonzero element of $\syz(p)_{2,3}$.  As
explained in the proof of the theorem, this is the special syzygy
whose multiples by polynomials in $z,w$ account for the question marks
in \eqref{mnpict}.

The corresponding curve $W(p)\cap Y  \subset W$ is parametrized by
\begin{align*}
\varphi(s,t) &= t^{2} \, p_{0} + s^{2}\, p_{1} - st \,p_{2}\\ 
      &= t^{2}\cdot x^{2}z + s^{2} \cdot y^{2}w - st (x^{2}w + 
           y^{2}z)\\
      &= (tx^{2} - sy^{2})(tz - sw).
\end{align*}
The third line shows that we have a curve in $W(p) \cap Y$, and
the first line makes it clear that it is a smooth conic.
\end{example}

\section{Hilbert Functions and Generators}
\label{genhilb}

\noindent {\bfseries \emph{Tools}:} Hilbert functions, functoriality

\subsection{Hilbert Functions}Theorem~\ref{th:main}
tells us that triples $(p_0,p_1,p_2) \in R_{2,1}^3$ with no common
zeros on $\P^1 \times \P^1$ come in two flavors, generic and
non-generic.  It is now easy to determine the Hilbert function of
$\syz(p)$.  

The Hilbert function of a finitely generated bigraded
$R$-module $M$ is defined by
\[
H_M(m,n) := \dim_k M_{m,n}.
\]
For example, the Hilbert function of $R = k[x,y,z,w]$ is
\[
H_R(m,n) = \begin{cases} (m+1)(n+1) & m,n \ge -1\\ 0 &
  \text{otherwise}.\end{cases}
\]
Here is the Hilbert function of the syzygy module.

\begin{proposition}
\label{hilbsyz} 
Let $\syz(p)$ be the syzygy module on $p_0,p_1,p_2$.
\begin{enumerate}
\item In the generic case, $\syz(p)$ has Hilbert function
\[
H_{\syz(p)}(m,n) = \begin{cases} 
3(m-3)(n-1) - (m-5)(n-2) & m \ge 5, n \ge 2\\
m-5 & m \ge 6, n = 1\\
2(n-2) & m = 3, n \ge 3\\
3(n-1) + n-2 & m = 4, n \ge 3\\
0 & \text{\rm otherwise}.
\end{cases}
\]
\item In the non-generic case, $\syz(p)$ has has Hilbert function
\[
H_{\syz(p)}(m,n) = \begin{cases} 
3(m-3)(n-1) - (m-5)(n-2) & m \ge 5, n \ge 2\\
m-5 & m \ge 6, n = 1\\
(m-1)(n-2) & m = 2,3, n \ge 3\\
3(n-1) + n-2 & m = 4, n \ge 3\\
0 & \text{\rm otherwise}.
\end{cases}
\]
\end{enumerate}
\end{proposition}

\begin{proof}
When $m \ge 5, n \ge 2$, Lemma~\ref{syzE2} and
Proposition~\ref{easypart} give the exact sequence
\[
0 \longrightarrow R_{m-6,n-3} \longrightarrow R_{m-4,n-2}^3
\longrightarrow \syz(p)_{m,n} \longrightarrow 0.
\]
This easily implies that $H_{\syz(p)}(m,n) = 3(m-3)(n-1) - (m-5)(n-2)$
for these $(m,n)$'s.

When $n = 1$, $m = 3$, or $m = 4$, the formulas for $H_{\syz(p)}(m,n)$
follow similarly from Lemma~\ref{syzE2} and
Proposition~\ref{easypart}.  When $m = 2$ and $n \ge 3$,
Theorem~\ref{th:main} tells us that $H_{\syz(p)}(2,n)$ equals $0$ in
the generic case and $n-2$ in the non-generic case.

The remaining cases are all $0$ by the zeros in \eqref{mnpict}.
\end{proof}

Using Proposition~\ref{hilbsyz} and the exact sequence
\[
0 \longrightarrow \syz(p) \longrightarrow R(-2,-1)^3 \longrightarrow I
\longrightarrow 0,
\]
it is easy to determine the Hilbert function of the ideal $I = \langle
p_0,p_1,p_2\rangle$.  Hence $I$ also has two Hilbert functions,
depending on whether $p_0,p_1,p_2$ are generic or non-generic.

\subsection{Generators of the Syzygy Module} 
We next determine degrees of the minimal generators of the syzygy
module in the generic and non-generic cases.

\begin{proposition}
\label{gensyz} 
Let $\syz(p)$ be the syzygy module on $p_0,p_1,p_2$.
\begin{enumerate}
\item In the generic case, $\syz(p)$ has six minimal generators: one
in degree $(6,1)$, three in degree $(4,2)$, and two in degree $(3,3)$.
\item In the non-generic case, $\syz(p)$ has five minimal generators:
one in degree $(6,1)$, three in degree $(4,2)$, and one in degree
$(2,3)$.
\end{enumerate}
\end{proposition}

\begin{proof}
First assume that $p_0,p_1,p_2$ are generic.  Then the question marks
in \eqref{mnpict} are all zeros.  Furthermore, we have the three
Koszul syzygies in degree $(4,2)$, which generate $\syz(p)$ in degree
$(m,n)$ for all bullets in \eqref{mnpict}.

It remains to consider the boxes.  When $n = 1$, we get the horizontal
line of boxes in \eqref{mnpict}.  We proved that $\syz(p)_{m,1}$ has
dimension $m-5$ for $m \ge 6$ using
\begin{align*}
\syz(p)_{m,1} = {}''\!E_2^{0,1}(m,n) &= H^1(m-6,-2)\\ &\simeq
H^0(\P^1,\Oc_{\P^1}(m-6)) \otimes H^1(\P^1,\Oc_{\P^1}(-2))\\ &\simeq
k[x,y]_{m-6}.
\end{align*}
The key point is that everything here---spectral sequences, the
K\"unneth formula, etc.---is functorial.  In particular, the above
isomorphisms are compatible with multiplication by homogeneous
polynomials in $k[x,y]$.  This gives the commutative diagram
\[
\begin{matrix}
\syz(p)_{6,1}\otimes k[x,y]_{m-6} &\simeq & k[x,y]_{0}\otimes
k[x,y]_{m-6}\\[2pt] \downarrow\ \ && \downarrow\ \ \ \,\\[2pt]
\syz(p)_{m,1} &\simeq & k[x,y]_{m-6}.
\end{matrix}
\]
where the vertical maps are multiplication.  Since the vertical map on
the right is an isomorphism, the same is true for the map on the
left.  This shows that the non-Koszul syzygy in degree $(6,1)$
generates $\syz(p)_{m,1}$ for $m \ge 6$.

A similar argument shows that the two non-Koszul syzygies in
degree $(3,3)$ generate $\syz(p)_{3,n}$ for $n \ge 3$.  For $m = 4$,
this argument also shows that 
\[
{}''\!E_2^{0,1}(4,3) \otimes k[z,w]_{n-3} \longrightarrow
{}''\!E_2^{0,1}(4,n)
\]
is surjective.  Recall from Proposition~\ref{easypart} that
${}''\!E_2^{0,1}(4,3)$ is one-dimensional, corresponding to one
non-Koszul syzygy in degree $(4,3)$.  Using Lemma~\ref{syzE2}, we see
that $\syz(p)_{4,n}$ is generated by the Koszul syzygies in degree
$(4,2)$ plus the non-Koszul syzygy in degree $(4,3)$.

It remains to show that the non-Koszul syzygy in degree $(4,3)$ can be
generated using the non-Koszul syzygies in degree $(3,3)$.  For this,
we first observe that 
\begin{equation}
\label{e2h1}
\begin{aligned}
{}''\!E_2^{0,1}(3,3) &\simeq H^1(\P^1,\Oc_{\P^1}(-3)) \otimes
H^0(\P^1,\Oc_{\P^1}(0))\\
{}''\!E_2^{0,1}(4,3) &\simeq H^1(\P^1,\Oc_{\P^1}(-2)) \otimes
H^0(\P^1,\Oc_{\P^1}(0)).
\end{aligned}
\end{equation}
Furthermore, since the duality \eqref{sd} is functorial, the map
\begin{equation}
\label{dualmap}
H^1(\P^1,\Oc_{\P^1}(-3)) \otimes k[x,y]_1 \longrightarrow
H^1(\P^1,\Oc_{\P^1}(-2)) 
\end{equation}
is dual to 
\[
H^0(\P^1,\Oc_{\P^1}(0)) \otimes k[x,y]_1 \longrightarrow
H^0(\P^1,\Oc_{\P^1}(1)).
\]
This map is $k[x,y]_0 \otimes k[x,y]_1 \to k[x,y]_1$, which is
obviously injective.  Hence the dual map \eqref{dualmap} is
surjective, which by the functorality of the isomorphisms \eqref{e2h1}
means that the non-Koszul syzygy in degree $(4,3)$ comes from the
syzygies in degree $(3,3)$.  This completes the proof in the generic
case. 

In the non-generic case, we get a non-Koszul syzygy in degree $(2,3)$.
The proof of Theorem~\ref{th:main} shows that it generates
$\syz(p)_{2,n}$ for $n \ge 3$.  Multiplying this syzygy by $x$ and $y$
gives two syzygies in degree $(3,3)$.  It is easy to see that these
are linearly independent and hence generate $\syz(p)_{3,3}$.  Then the
argument in the generic case shows that this syzygy, together with the
Koszul syzygies in degree $(4,2)$ and the non-Koszul syzygy in degree
$(6,1)$, generate $\syz(p)$.
\end{proof}

\subsection{The Free Resolution of $I$}  If we combine the generators
of the syzygy module described in Proposition~\ref{gensyz} with the
exact sequence
\[
0 \longrightarrow \syz(p) \longrightarrow R(-2,-1)^3 \longrightarrow I
\longrightarrow 0,
\]
then we can start to build a minimal free resolution of the  ideal $I =
\langle p_0,p_1,p_2\rangle$.  In the generic case, we get
\begin{equation}
\label{startg}
\cdots \longrightarrow
\begin{array}{c}
 R(-6,-1)\\
 \oplus \\
R(-4,-2)^3\\
\oplus\\
R(-3,-3)^2 \\
\end{array}
 \longrightarrow
R(-2,-1) ^3
\longrightarrow I \longrightarrow 0,
\end{equation}
while in the non-generic case, the resolution begins
\begin{equation}
\label{startng}
\cdots
 \longrightarrow
\begin{array}{c}
 R(-6,-1)\\
 \oplus \\
R(-4,-2)^3\\
\oplus\\
R(-2,-3) \\
\end{array}
 \longrightarrow
R(-2,-1) ^3
\longrightarrow I \longrightarrow 0.
\end{equation}
In Section~\ref{shape} we will show how to complete these partial
resolutions.  

Recall that once we have a minimal free resolution, it is
straightforward to determine the Hilbert function.  Since we already
know the Hilbert function (Proposition~\ref{hilbsyz}), one might hope
to reconstruct the minimal free resolution from the Hilbert function.
Unfortunately, this is not possible in general, since different
minimal free resolutions can have the same Hilbert function.

\begin{example}
Let $S = k[x]$, graded in the usual way.  Then consider the
graded $S$-modules $M_1$ and $M_2$ defined by the exact sequences
\[
0 \longrightarrow S(-2) \xrightarrow{\left[{\displaystyle x^2}\right]} S
\longrightarrow M_1 \longrightarrow 0
\]
and 
\[
0 \longrightarrow 
\begin{array}{c}
 S(-2)\\
 \oplus \\
 S(-1)
\end{array}
\xrightarrow{\left[\begin{matrix} x & 0\\ 0 & x
\end{matrix}\right]}
\begin{array}{c}
 S(-1)\\
 \oplus \\
 S
\end{array}
\longrightarrow M_2 \longrightarrow 0.
\]
These sequences give the minimal free resolutions of $M_1$ and $M_2$.
One easily sees that $M_1$ and $M_2$ have the same Hilbert
function, even though the free resolutions are very different.
\end{example}

We will see in Section~\ref{shape} that in our situation, there is a
one-to-one correspondence between Hilbert functions and minimal free
resolutions.

\section{Explicit Non-Koszul Syzygies}
\label{explicit}

\noindent {\bfseries \emph{Tools}:} determinants, resultants

\subsection{Formulas for Syzygies} Koszul syzygies have the nice
property that they are given by explicit formulas; namely, they are
generated by the columns
of the $3\times3$ matrix in \eqref{koszul}.  The goal of this section
is to give explicit formulas for the non-Koszul syzygies.

Suppose that for  $i = 0,1,2$ we have polynomials
\begin{equation}
\label{pifor} 
p_{i} = {a_i}\,x^{2}z + {b_i}\,xyz + {c_i}\,y^2z + {d_i}\,x^{2}w +
{e_i}\,xyw + {f_i} \, y^2w
\end{equation}
with no common zeros on $\mathbb{P}^{1} \times \mathbb{P}^{1}$.  We
will describe a simple method for constructing one explicit nonzero
syzygy in degree $(6,1)$ and two linearly independent syzygies in
degree $(3,3)$.  

We start with a straightforward lemma.

\begin{lemma}  
\label{lemma:gh}
Suppose that there exist bihomogeneous polynomials $g,h$ such that
\begin{equation} 
\label{eq:gh}
p_{i}   =  A_i \, g  +  B_i \, h, \quad i =0,1,2,
\end{equation}
and let $C_{ij} = A_i B_j - A_j B_i$.  Then $C := (C_{12}, C_{20},
C_{01})$ is a syzygy on $p_{0}, p_{1}, p_{2}$.
\end{lemma}

\begin{proof}
The determinant
\[  
\det\left[\begin{matrix}
p_{0} & A_0 & B_0 \\
p_{1} & A_1 & B_1 \\
p_{2} & A_2 & B_2 
\end{matrix}\right]
\]
vanishes, since by our assumption \eqref{eq:gh}, its first column
equals $g$ times its second column plus $h$ times the third.  The
Laplace expansion of this determinant along the first column gives
the identity
\[ 
p_{0}\, (A_1 B_2 - A_2 B_1) + p_{1} \, (A_2 B_0 - A_0 B_2) + p_{2} \,
( A_0 B_1 - A_1 B_0) = 0.
\]
This proves that $C$ is a syzygy on $p_{0}, p_{1}, p_{2}$.
\end{proof}

This lemma can be used to construct Koszul syzygies as follows.

\begin{example}
\label{rmk:Koszul}
If we let $g = 1$, $h = 0$, and $A_i = p_{i}$ in \eqref{eq:gh}, then
Lemma~\ref{lemma:gh} gives a syzygy for any choice $B_0$, $B_1$, $B_2$
of bihomogeneous polynomials of the same degree.  One easily sees
that the resulting syzygy is 
\[
B_0\,(0,-p_2,p_1) + B_1\,(-p_2,0,p_0) + B_2\,(p_1,-p_0,0).
\]
This is an arbitrary Koszul syzygy.  In other words, letting
$(B_0,B_1,B_2) = (1,0,0)$, $(0,1,0)$, and $(0,0,1)$ gives generators
of the submodule of Koszul syzygies.  
\end{example}

We now construct a syzygy of degree $(6,1)$.  Recall that this will
involve polynomials of degree $(4,0)$ since $(6,1) = (4,0) + (2,1)$.
Using \eqref{pifor}, we write each $p_i$ in the form
\[ 
p_{i} =  A_i(x,y) \, z  +  B_i(x,y) \, w, \quad  i=0,1,2, 
\]
with 
\[ 
A_i (x, y )  =  a_i \, x^2  +  b_i \, x y  + c_i \, y ^2, \quad 
B_i (x, y )  =  d_i \, x^2  +  e_i \, x y  + f_i \, y ^2
\]
of degree $(2,0)$.  Thus the syzygy $C$ from Lemma~\eqref{lemma:gh}
involves polynomials of degree $(4,0)$.

\begin{proposition}
\label{prop:deg61}
Let $A_i$ and $B_i$ be as above   Then the syzygy $C$ from
Lemma~\ref{lemma:gh} is a non-Koszul syzygy in $\syz(p)_{6,1}$.
\end{proposition}

\begin{proof} 
Taking $g=z, h=w$ in the previous lemma, we know that $C$ is indeed a
syzygy.  Note that all $C_{ij}$ have degree $(4,0)$ and so they cannot
be a multiple of any $p_{i}$.  So, to prove that the syzygy is
non-Koszul we only need to show that $C$ is nonzero.  If this were so,
then the matrix
\[
\left[\begin{matrix}
A_0 & A_1 & A_2\\
B_0 & B_1 & B_2
\end{matrix}\right]
\]
would have rank at most $1$.  Hence there would be $(\alpha, \beta)
\in \A^2 \setminus \{(0,0)\}$ such that $\alpha (A_0, A_1, A_2) +
\beta (B_0, B_1, B_2) = 0$.  Then the three polynomials $p_{0}, p_{1},
p_{2}$ would vanish at any point of the form $(x,y,\alpha,\beta)$, a
contradiction.
\end{proof}

We next construct two syzygies of degree $(3,3)$.  These will involve
polynomials of degree $(1,2)$ since $(3,3) = (1,2) + (2,1)$.  Write
the polynomials $p_0,p_1,p_2$ from \eqref{pifor} in the form
\begin{equation}
\label{cdedef}
p_{i} = C_i(z,w)\, x^2 + D_i(z,w) \, x y + E_i(z,w) \, y^2, \quad
i=0,1,2,
\end{equation}
where $C_i,D_i,E_i$ are homogeneous linear forms in $z,w$, i.e.,
bihomogeneous polynomials of degree $(0,1)$.  Explicitly, $C_i(z,w) =
a_i \, z + d_i \, w$, $D_i(z,w) = b_i \, z + e_i \, w$, and $E_i(z,w)
= c_i \, z + f_i \, w$.  If we take $g = x, h = y^2$ in
Lemma~\ref{lemma:gh} and write
\[ 
p_{i} = \left(C_i(z,w) \, x + D_i(z,w) \, y\right) \, x + E_i(z,w)
\, y^2, \quad i=0,1,2,
\]
then we get a syzygy $C^{(1)}$.  On the other hand, if we take $g =
x^2, h =y$ and write
\[ 
p_{i} = C_i(z,w) \, x^2 + \left( D_i(z,w) \, x + E_i(z,w) \, y\right) \, y,
\quad i=0,1,2,
\]
then Lemma \ref{lemma:gh} gives a second syzygy $C^{(2)}$.   Note that
$C^{(1)}$ and $C^{(2)}$ lie in $\syz(p)_{3,3}$.

\begin{proposition} 
\label{prop:deg33}
The non-Koszul syzygies in degree $(3,3)$ are spanned by the syzygies
$C^{(1)}$ and $C^{(2)}$ constructed above.
\end{proposition}

\begin{proof}
Since $\dim_k \syz(p)_{3,3} = 2$ and there are no Koszul syzygies in
degree $(3,3)$, it is enough to show that $C^{(1)}$ and $C^{(2)}$ are
linearly independent.

We first observe that $C^{(1)} \ne 0$.  To see why, note that $C^{(1)}
= 0$ implies that
\begin{align*}
(C_1x+D_1y)E_2 - (C_2x+D_2y)E_1 &= 0\\
(C_2x+D_2y)E_0 - (C_0x+D_0y)E_2 &= 0\\
(C_0x+D_0y)E_1 - (C_1x+D_1y)E_0 &= 0.
\end{align*}
The first equation can be rewritten as 
\[
(C_1E_2 - C_2E_1)\,x - (D_1E_2-D_2E_1)\,y = 0,
\]
and similarly for the other equations.  Hence these equations imply
the vanishing of certain $2\times2$ minors of
\begin{equation}
\label{cdemat}
\left[\begin{matrix}
C_0 & C_1 & C_2 \\
D_0 & D_1 & D_2\\
E_0 & E_1 & E_2 
\end{matrix}\right].
\end{equation}
It follows easily that \eqref{cdemat} does not have maximal rank,
which by \eqref{cdedef} implies that $p_{0}, p_{1}, p_{2}$ are
linearly dependent.  Yet we observed in Section~\ref{QM} that $p_{0},
p_{1}, p_{2}$ are linearly independent since they have no common zeros
on $\P^1\times\P^1$.  Thus $C^{(1)} \ne 0$.

If $C^{(1)}$ and $C^{(2)}$ were linearly dependent, then $C^{(2)}$
would be a constant times $C^{(1)}$.  Since $C^{(1)}$ is formed from
the $2\times2$ minors of
\[
\left[\begin{matrix}
C_0 \,x + D_0 \,y & C_1 \,x + D_1\, y & C_2\, x + D_2\, y \\
E_0 & E_1 & E_2
\end{matrix}\right], 
\]
it lies in the kernel of this matrix, and similarly $C^{(2)}$ lies in
the kernel of
\[
\left[\begin{matrix}
C_0 & C_1 & C_2\\
D_0 \,x+ E_0 \,y & D_1 \,x + E_1\, y & D_2\, x + E_2\, y 
\end{matrix}\right].
\]
If $C^{(2)}$ were a multiple of $C^{(1)}$, then it is straightforward
to show that $C^{(1)}$ would lie in the kernel of \eqref{cdemat}.
Since $C^{(1)} \ne 0$, we conclude that \eqref{cdemat} would not have
maximal rank.  As above, this is impossible, and linear independence
follows.
\end{proof}

Propositions~\ref{prop:deg61} and~\ref{prop:deg33} give explicit
formulas for one non-Koszul syzygy of degree $(6,1)$ and two linearly
independent non-Koszul syzygies of degree $(3,3)$.  When combined with
the three Koszul syzygies of degree $(4,2)$, Proposition~\ref{gensyz}
shows that we have explicit formulas for the generators of the syzygy
module $\syz(p)$ in the generic case.  In Section~\ref{shape}, we will
use the internal structure of these explicit formulas to determine the
free resolution of $I = \langle p_0,p_1,p_2\rangle$ in the generic
case.

\subsection{The Non-Generic Case}
As we have already seen, non-Koszul syzygies in degree $(2,3)$ occur
only when $p_0,p_1,p_2$ are non-generic.  They correspond to the case
when the linear span of our polynomials intersects the variety of
$R_{2,0}\times R_{0,1}$ factorizable polynomials in a smooth conic.

One way to proceed is to compute the syzygies $C^{(1)}, C^{(2)}$ of
degree $(3,3)$.  The analysis of Section~\ref{genhilb} shows that
these are multiples of a syzygy of degree $(2,3)$ by linear forms in
$x,y$.  Then we get the desired syzygy of degree $(2,3)$ by taking the
GCD of $C^{(1)}$ and $C^{(2)}$.

\begin{example} We revisit Example \ref{ex:special}, where
$p_{0} = x^{2}z$, $p_{1} = y^{2}w$ and $p_{2} = x^{2}w + 
y^{2}z$.  One computes that
\begin{align*}
C^{(1)} &= (-xw^2,-xz^2,xzw)\\
C^{(2)} &= (-yw^2,-yz^2,yzw).
\end{align*}
These are multiples (by $-x$ and $-y$ respectively) of
the syzygy $(w^2,z^2,-zw)$ described in Example \ref{ex:special}.
\end{example}

In general, we can use this method to construct the non-Koszul syzygy
in degree $(2,3)$ when $W(p) \cap Y$ is a smooth conic.  We first state
an easy consequence of Theorem~\ref{th:main} that we will also need
in Section \ref{shape}.

\begin{lemma} \label{lem:cigens}
In the non-generic case, we may assume after a change of basis in
$W(p) = \P(\mathrm{Span}(p_0,p_1,p_2))$ that $p_0, p_1, p_2$ are of
the form $p_i = g_i l_i, \, i=0,1,2,$ where $g_i \in R_{(2,0)}$ and
$l_i \in R_{(0,1)}$.  Furthermore, we may assume that $g_2 = g_0 +
g_1$, $l_0 =z$, $l_1 =w$, and $l_2 = z+w$.  Moreover, $\langle g_0,g_1
\rangle$ and $\langle p_0, p_1\rangle$ are complete intersections.
 \end{lemma}

\begin{proof}
As in the proof of Theorem~\ref{th:main}, we can parametrize $W(p) \cap Y$ by
\[
\phi(s,t) = (r_0(s,t)x^2 + r_1(s,t)xy + r_2(s,t)y^2)(tz - sw)
\]
where $r_i(s,t)$ is linear in $s,t$ (see \eqref{vpfor}).  
Note that any three points on a smooth conic define the same $\P^2$,
since $n+1$ points on a rational normal curve of degree $n$ are in
linearly general position by \cite[p.~10]{harris}.
Hence we can replace
our original $p_0,p_1,p_2$ with
\begin{equation}
\label{specialp}
p_0 = \phi(0,1), p_1 = \phi(-1,0),\ p_2 = \phi(-1,1),
\end{equation}
i.e., $p_i = g_i l_i$ with $l_0 =z, l_1 =w, g_0 =
r_0(0,1)x^2+r_1(0,1)xy+r_2(0,1)y^2, g_1=
r_0(-1,0)x^2+r_1(-1,0)xy+r_2(-1,0)y^2$ and $l_2 = z+w, g_2 = g_0+g_1,
$ by the linearity of the $r_i$.  If $ \langle g_0, g_1 \rangle$ is
not a complete intersection, they share a zero $(x_0,y_0) \in \P^1$,
from which $(x_0, y_0, z,w)$ would be a common zero of $p_0, p_1, p_2$
for any value of $(z,w)$, a contradiction. It easily follows that
$\langle p_0, p_1 \rangle$ is also a complete intersection.
\end{proof}

Changing to a different basis of $W(p) = \P(\mathrm{Span}(p_0,p_1,p_2))$ 
induces an isomorphism of the corresponding syzygy modules. So, we assume
that our non-generic polynomials $p_0,p_1,p_2$ are as in the statement of 
Lemma \ref{lem:cigens} and we keep the notation of the proof.
One can compute that in this case the syzygies $C^{(1)}$ and $C^{(2)}$ of
Proposition~\ref{prop:deg33} are given by
\begin{align*}
C^{(1)} &= -(m_{13}\, x + m_{23}\, y)(w(z+w), z(z+w), -zw)\\
C^{(2)} &= -(m_{12}\, x + m_{13}\, y)(w(z+w), z(z+w), -zw),
\end{align*}
where $m_{ij}$ is the $2\times2$ minor formed using columns $i$ and $j$ of
the rank $2$ matrix
\[
\left[\begin{matrix} 
r_0(0,1) & r_1(0,1) & r_2(0,1)\\
r_0(-1,0) & r_1(-1,0) & r_2(-1,0)
\end{matrix}\right].
\]
Here, $C^{(1)}$ and $C^{(2)}$ are common multiples of $(w(z+w),
z(z+w), -zw)$, which is easily seen to be a syzygy on \eqref{specialp}
since $r_i(-1,1) = r_i(0,1) + r_i(-1,0)$.  This gives the desired
element of $\syz(p)_{2,3}$.

\subsection{Resultants} \label{subs:resultants}
We end this section with a ``computational''
translation of our hypothesis that $p_{0}, p_{1}, p_{2}$ have no
common zeros on $\mathbb{P}^{1} \times \mathbb{P}^{1}$.  This
condition is precisely described by the non-vanishing of the
multihomogeneous unmixed resultant
$\mathrm{Res}_{(2,1)}(p_0,p_1,p_2)$.  As described in \cite[Ch.\
7]{clo}, we can regard this resultant as an irreducible polynomial
$\mathrm{Res}_{(2,1)} \in \Z[a_0,\dots,f_2]$, where the variables
$a_0,\dots,f_2$ represent the $18$ coefficients of $p_{0}, p_{1},
p_{2}$ in \eqref{pifor}.

The polynomial $\mathrm{Res}_{(2,1)}$ can be expressed most
efficiently as a determinant.  Of the known ways of doing this, the
smallest can be obtained using the methods in \cite{de}.  Consider the
dehomogenizations $P_{i}(x,z) := p_i (x,1,z,1)$.  Thus
\[
P_{i} =   {a_i}\,x^{2}z +  {b_i}\,x z + 
 {c_i}\,z +  {d_i}\,x^{2} +  {e_i}\,x + 
 {f_i}  \, \in k[x,z], \quad i=0,1,2.
\]
Let $X,Z$ be two new variables, and denote by $B$ the matrix
\[
B  = 
\left[\begin{matrix}
P_0(x,z) & P_1(x,z) & P_2(x,z)\\
P_0(X,z) & P_1(X,z) & P_2(X,z)\\
P_0(X,Z) & P_1(X,Z) & P_2(X,Z)
\end{matrix}\right].
\]
One can show that the expression 
\[
b = \frac 1 {(x-X) (z-Z)} \det(B)
\] 
is in fact a polynomial in $x$, $z$, $X$, and $a_0,\dots,f_2$.  It is
also easy to see that $b$ has degree $3$ in $X$ and degree $1$ in $x$
and $z$.  Thus $b$ has an expansion
\[
b(x,z,X) = \sum_{j=0}^3 b_{0j}\,  X^j +  b_{1j}\,  x X^j +  
b_{2j} \, z X^j +  
b_{3j} \, x z X^j,
\]
where $b_{ij} \in \Z[a_0,\dots,f_2]$ for $0 \le i,j \le 3$.  The
results of \cite{de} imply that ${\rm Res}_{(2,1)} = \det(b_{ij})$.
This compact expression is a polynomial of degree $12$ in the $18$
coefficients of $p_0,p_1,p_2$, whose expansion as a sum of monomials
in $a_0,\dots,f_2$ has $20,791$ terms.

\section{The Shape of the Minimal Free Resolutions}
\label{shape}

\noindent {\bfseries \emph{Tools}:} computations, mapping
cones, Hilbert-Burch theorem, liason

\subsection{Computations} 
In Section 5.3, we determined the first syzygies of an ideal generated
by three elements of $R_{(2,1)}$ with no common zeros on
$\P^1\times\P^1$; the syzygies differ in the generic and non-generic
cases. It is natural to ask about higher syzygies, and more generally
about the minimal free resolution.  Resolutions can be calculated by a
computer, which we illustrate using Macaulay~2.

\begin{example}
\label{specialex}
For multigraded rings, Macaulay 2 requires weights which have a
nonzero first entry.  So we will alter the grading of $R$ a bit,
creating a ring where $x,y$ have degree $(1,0)$ and $z,w$ have degree
$(1,1)$.  Hence $p_0,p_1,p_2$ have degree $(3,1)$ in this grading.  

We return again to Example \ref{ex:special}:

\medskip

\begin{small}
\renewcommand{\baselinestretch}{.8}
\begin{verbatim}
i1 : R = ZZ/31991[x,y,z,w,Degrees => {{1,0},{1,0},{1,1},{1,1}}];

i2 : I = ideal (x^2*z,y^2*w,x^2*w+y^2*z)

             2    2    2     2
o2 = ideal (x z, y w, y z + x w)

o2 : Ideal of R

i3 : res I  

-- The res command produces a minimal free resolution.

      1      3      5      4      1
o3 = R  <-- R  <-- R  <-- R  <-- R  <-- 0

     0      1      2      3      4      5

-- To see all the differentials in a complex, append the suffix .dd: 

i4 : o3.dd

          1                           3
o4 = 0 : R  <----------------------- R  : 1
               | x2z y2z+x2w y2w |

          3                                                  5
     1 : R  <---------------------------------------------- R  : 2
               {3, 1} | -w2 -y2z-x2w -y2w 0        -y4  |
               {3, 1} | zw  x2z      0    y2w      x2y2 |
               {3, 1} | -z2 0        x2z  -y2z-x2w -x4  |

          5                                 4
     2 : R  <----------------------------- R  : 3
               {5, 3} | x2 -y2 0   0   |
               {6, 2} | -w 0   -y2 0   |
               {6, 2} | z  w   x2  -y2 |
               {6, 2} | 0  z   0   -x2 |
               {7, 1} | 0  0   z   w   |

          4                     1
     3 : R  <----------------- R  : 4
               {7, 3} | y2 |
               {7, 3} | x2 |
               {8, 2} | -w |
               {8, 2} | z  |

          1
     4 : R  <----- 0 : 5
               0
\end{verbatim}
\end{small}

\medskip

In this example, $p_0= x^2z$,  $p_1= y^2z+x^2w$, 
$p_2= y^2w$, and the second, third and fourth columns of the matrix
of the second differential correspond to the generators of the Koszul syzygies
$k_{01}=(-p_1,p_0,0), k_{02}=(-p_2, 0, p_0)$ and $k_{21}=(0, p_2, -p_1)$.
Also, calling $C_1, \dots, C_4$ the columns of the matrix of the
third differential (which generate the second syzygies), note that
$y^2 C_1 + x^2 C_2 = w C_3 - z C_4 = (0, -p_2, p_1, p_0, 0)$, yielding
the second Koszul syzygy $(-p_2, p_1, p_0)$ 
of the first Koszul syzygies $ -p_2 k_{01} + p_1 k_{02} + p_0 k_{21} =0$.
\end{example}

\begin{example}
\label{genericex}
We now consider the generic case.  We use the same setup as for
Example~\ref{specialex}, except that now we ask Macaulay 2 to resolve
the ideal generated by 3 random elements of degree (3,1).  

We obtain the following resolution:
\medskip
\begin{small}
\begin{verbatim}
i5 : res ideal random(R^{3:{3,1}},R^1)

      1      3      6      5      1
o5 = R  <-- R  <-- R  <-- R  <-- R  <-- 0

     0      1      2      3      4      5

o5 : ChainComplex

-- We do not display the first two differentials, which are dense
-- and big.  d_3 and d_4 are more compact, hinting the resolution
-- may have a nice structure.  For compactness, some coefficients 
-- of d_3 are written as *.

i6 : o5.dd_3

o6 = {6, 2} | -z+*w   *w      *w       *x2+*xy+*y2  *xy+*y2     |
     {6, 2} | *w      -z+*w   -z+*w    *xy+*y2      *x2+*xy+*y2 |
     {6, 2} | *w      *w      -z+*w    *x2+*xy+*y2  *xy+*y2     |
     {6, 3} | x       y       0        0            0           |
     {6, 3} | 10287y  -6711y  x+8060y  0            0           |
     {7, 1} | 0       0       0        -z+13959w    -w          |

             6       5
o6 : Matrix R  <--- R


i7 : o5.dd_4

o7 = {7, 3} | 10780xy-756y2    |
     {7, 3} | -10780x2+756xy   |
     {7, 3} | -12929xy-13054y2 |
     {8, 2} | w                |
     {8, 2} | -z+13959w        |
\end{verbatim}
\end{small}
\renewcommand{\baselinestretch}{1.0}
\end{example}

\medskip

In Proposition~\ref{gensyz}, we learned that there are only two
possible types of behavior for the first syzygies.  In general, this
tells us nothing about higher syzygies.  However, the zeros in the
third differentials in Examples~\ref{specialex} and~\ref{genericex}
suggest that there may be some interesting structure in the
resolutions.  The intuition is that having a zero entry means that
some syzygy involves only part of the previous syzygy module, so there
may be a way to split things up.  This is indeed the case, and we now
discuss the tools we shall need.

\subsection{The mapping cone}
Let $R$ be a ring, $f \in R$ and $L$ an ideal of $R$. 
There is an exact sequence 
\[
0 \longrightarrow (L+\langle f \rangle)/L \longrightarrow  
R/L \longrightarrow R/(L+\langle f \rangle) \longrightarrow 0.
\]
The map $1 \mapsto f$ gives a surjection from $R$ to $(L+\langle f
\rangle)/L$; the kernel of the map is simply $\{g \in R \mid fg \in
L\} = (L:f)$. Notice that if $R$ is graded, $L$ is homogeneous, and $f
\in R_a$, then we actually obtain a graded exact sequence:
\[
0 \longrightarrow R(-a)/(L:f) \stackrel{[f]}{\longrightarrow}  
R/L \longrightarrow R/(L+\langle f \rangle) \longrightarrow 0.
\]
Suppose we have graded free resolutions $F_\bullet$ of $R(-a)/(L:f)$
and $G_\bullet$ of $R/L$ with $F_0 = R(-a)$ and $G_0 = R$.  Then we
have a commutative diagram:

\vspace*{.89in}

\begin{equation}
\label{beginnings}
\end{equation}

\vspace*{-1.55in}

\[
\xymatrix{
 &&&& 0 \ar[d] \\
\ar[r]^{d_3} &F_2 \ar[r]^{d_2} & F_1 \ar[r]^{d_1} & R(-a) \ar[r] \ar[d]^{[f]} &
 R(-a)/(L:f) \ar[d]^{[f]} \ar[r]& 0\\
\ar[r]^{\delta_3}& G_2 \ar[r]^{\delta_2} & G_1 \ar[r]^{\delta_1} &
  R \ar[r] \ar[dr] & R/L \ar[d]\ar[r] &0\\
&&&&R/(L+\langle f \rangle) \ar[d]\\
&&&& 0 \\
}
\]
\vskip1pt
\noindent The rightmost column is exact by construction.  Furthermore,
an easy diagram chase shows that the diagonal map fits into an exact
sequence
\[
R(-a) \oplus G_1 \stackrel{[f] \oplus \delta_1}{\longrightarrow} R
\longrightarrow R/(L+\langle f \rangle) \longrightarrow 0.
\]
So we have the beginning of a free resolution for $R/(L+\langle f
\rangle)$. 

\begin{lemma}
\label{lem:FGmap}
In the setting above, there is a map of complexes from $F_\bullet$ to
$G_\bullet$ such that $H_0(F_\bullet) \to H_0(G_\bullet)$ is the map
\[
R(-a)/(L:f) \stackrel{[f]}{\longrightarrow} R/L.
\]
\end{lemma}

\begin{proof}
This is a standard result in homological algebra---see \cite[Prop.\
A3.13]{e}, for example.  You should review the argument to see that if
the graded free resolutions $F_\bullet$ and $G_\bullet$ are given
explicitly, then so is the map between them.
\end{proof}

Now that we have a map $\psi_\bullet\colon F_\bullet \to G_\bullet$,
the process of constructing a resolution for $R/(L+\langle f \rangle)$
is pretty simple.  We already have the first two steps in the exact
sequence
\[
R(-a) \oplus G_1 \stackrel{[f] \oplus \delta_1}{\longrightarrow} R
\longrightarrow R/(L+\langle f \rangle) \longrightarrow 0.
\]
Let $M_0 := G_0 = R$, and for all $i \ge 1$, define modules $M_i :=
F_{i-1} \oplus G_i$, where $F_0 = R(-a)$.  Now define maps $M_i \stackrel
{\partial_i}{\longrightarrow} M_{i-1}$ via $\partial_1 := [f] \oplus
\delta_1$ for $i = 1$ and 
\[
\partial_i := \left[\begin{matrix} d_{i-1} & 0\\ (-1)^{i-1} \psi_{i-1}
    & \delta_i 
\end{matrix}\right]
\]
for $i \ge 2$.  A straightforward computation shows that $\partial_i
\circ \partial_{i-1} = 0$, so that $(M_\bullet, \partial_\bullet)$ is
a complex.  In this notation, the first two steps of the resolution of
$R/(L+\langle f \rangle)$ can be written as
\[
M_1 \longrightarrow M_0 \longrightarrow R/(L+\langle f \rangle)
\longrightarrow 0.
\]

\begin{lemma}\label{lem:mapcone}
The complex $(M_\bullet, \partial_\bullet)$ is resolution of
$R/(L+\langle f\rangle)$.
\end{lemma}

\begin{proof}
This can be proved by an explicit argument using the formula for
$\partial_i$ and the fact that $F_\bullet$ and $G_\bullet$ are
resolutions of $R(-a)/(L:f)$ and $R/L$ respectively.  Alternatively,
it is easy to see that there is a short exact sequence of complexes
\[
0 \longrightarrow G_\bullet \longrightarrow M_\bullet \longrightarrow
F_\bullet[-1] \longrightarrow 0,
\]
where $(F_\bullet[-1])_i = F_{i-1}$.  Using the corresponding
long exact sequence in homology (see \cite[Prop.\ A3.15]{e} or
\cite[Thm.\ 8.1.4]{s}) easily gives the desired result.  It is a good
exercise to work out the details of both proofs.
\end{proof}

The complex $(M_\bullet, \partial_\bullet)$ is called the
\emph{mapping cone}.  While it always gives a resolution, it need not
be minimal; in particular, not all the entries of $\partial_i$ need to
have positive degree. We illustrate the mapping cone with a pair of
examples.

\begin{example}
The minimal free resolution of the quotient of a ring by an ideal
generated by a regular sequence $f_1,\ldots, f_n$ is the Koszul
complex of the $f_i$.  This can be proved by induction using the
mapping cone construction since the definition of a regular sequence
implies that $(\langle f_1,\ldots, f_k\rangle :f_{k+1}) = \langle
f_1,\ldots, f_k\rangle$.  See \cite[\S 3.3]{s} for more details.
\end{example}

\begin{example}
\label{hbex}
In $k[x,y]$, consider the ideal $L = \langle x^2, y^2 \rangle$ and the
element $f = xy$.  It is easy to see that $(L:xy)=(x,y)$, so that
$F_\bullet$ is a Koszul complex on $x,y$, and $G_\bullet$ is a Koszul
complex on $x^2,y^2$.  The resulting mapping cone resolution of
$R/(L+\langle f\rangle) = R/\langle x^2,y^2,xy\rangle$ is therefore:
\[ 
0 \rightarrow R(-4) \rightarrow R(-4)\oplus R^2(-3) 
\rightarrow R^3(-2) \rightarrow R \rightarrow
R/\langle x^2,y^2,xy\rangle \rightarrow 0.
\]
However, the actual minimal resolution is:
\[ 
0 \longrightarrow  R^2(-3) 
\longrightarrow R^3(-2) \longrightarrow R \longrightarrow
R/\langle x^2,y^2,xy\rangle \longrightarrow 0,
\]
where the maps are given by:
\begin{small}
\renewcommand{\baselinestretch}{.8}
\begin{verbatim}

          1                    3
o3 = 0 : R  <---------------- R  : 1
               | x2 xy y2 |

          3                     2
     1 : R  <----------------- R  : 2
               {2} | -y 0  |
               {2} | x  -y |
               {2} | 0  x  |

          2
     2 : R  <----- 0 : 3
               0
\end{verbatim}
\end{small}
\renewcommand{\baselinestretch}{1.0}
\end{example}

In Example~\ref{hbex}, it is interesting to note that in the minimal
resolution, the two-by-two minors of the matrix of first syzygies
actually generate the ideal.  This is no accident---it is guaranteed
by the following classic result.

\begin{theorem}[The Hilbert-Burch theorem] Let
\label{hbst} 
\[
\mathcal{F}: 0 \longrightarrow F_2 \stackrel{d_2}{\longrightarrow} F_1  
\stackrel{d_1}{\longrightarrow} R \longrightarrow R/L \longrightarrow
0
\]
be a complex where $F_1$ and $F_2$ are free.  Then:
\begin{enumerate}
\item If $\mathcal{F}$ is exact and $F_1 \simeq R^n$, then $F_2 \simeq
R^{n-1}$ and there exists a nonzerodivisor $a$ such that $L = a
I_{n-1}(d_2)$, where $I_{n-1}(d_2)$ is the Fitting ideal generated by the
$(n-1)\times(n-1)$ minors of $d_2$.  Furthermore, the $i^{th}$ entry
of $d_1$ is $(-1)^ia$ times the minor obtained from $d_2$ by deleting
the $i^{th}$ row and the ideal $L$ has depth exactly two.
\item Conversely, given any $(n-1)\times n$ matrix $d_2$ with 
depth $I_{n-1}(d_2) \ge 2$, and a nonzerodivisor $a$, the map $d_1$ 
obtained as in part $(1)$ makes $\mathcal{F}$ into a free resolution of
$R/L$, with $L = a I_{n-1}(d_2)$.
\end{enumerate}
\end{theorem}

For a proof, we refer the reader to \cite[Theorem 20.15]{e}.  The
Hilbert-Burch theorem will prove to be useful below.

\subsection{The Non-Generic Case}
We now return to the problem of describing the free resolution of the
ideal $I = \langle p_0,p_1,p_2\rangle$.  We begin with the non-generic
case.  By Lemma~\ref{lem:cigens}, we may assume that
\[
p_0 = g_0 z,\ p_1 = g_1 w,\ p_2 = g_2 l_2 = (g_0 + g_1)(z+w),
\]
where $g_0,g_1 \in R_{2,0}$ and $C = \langle g_0z,g_1w\rangle$ is
a complete intersection.  We shall obtain the entire free resolution
of $I$ by taking two mapping cone resolutions.  Consider the short
exact sequence
\[
0 \longrightarrow R(-2,-1)/(C : g_2l_2) \stackrel{[g_2l_2]
}{\longrightarrow} R/C \longrightarrow R/I \longrightarrow 0.
\] 
Since $C$ is a complete intersection, the free resolution of $R/C$ is
a Koszul complex.  To obtain a mapping cone resolution for $R/I$
(which is not, in general, a minimal resolution), we have to determine
the ideal quotient $J = (C : g_2l_2)$.

The key to determining the ideal quotient is to use the fact that for
any ideal $L = \langle f_1,\ldots, f_k \rangle$ and element $f$, an
element $a$ is in $(L:f)$ if and only if there exist $a_1, \dots, a_k$
such that $a f = \sum_{i=1}^k a_if_i$, i.e. if and only if $(a_1,
\dots, a_k,-a)$ is a syzygy on $(f_1,\ldots,f_k,f)$.

In our situation, Proposition~\ref{gensyz} describes the degrees of
the minimal generators of the syzygy module.  The generators of
degrees $(4,2)$ are the Koszul syzgyies
\[
(p_1,p_0,0),\ (p_2,0-p_0),\ (0,p_2,-p_1),
\]
and by Proposition~\ref{prop:deg61}, the $2\times2$ minors of 
\[
\left[\begin{matrix} g_0 & 0 & g_2\\ 0 & g_1 & g_2\end{matrix}\right]
\]
give the generator $(-g_1g_2,-g_0g_2,g_0g_1)$ in degree $(6,1)$.
Also, we saw in Section~6.2 that in the situation of
Lemma~\ref{lem:cigens}, the generator in degree $(2,3)$ is the syzygy
$(w(z+w),z(z+w),-zw)$.  By the previous paragraph, it follows that $J
= (I : g_2l_2)$ is given by
\[
J = \langle zw, g_0z, g_1w, g_0g_1 \rangle. 
\]

We need to find the resolution of this ideal, and we do it with (yet
another!)\ mapping cone.  First, we'll let $K$ denote the ideal
$\langle zw, g_0z, g_1w \rangle$, so $J = K + \langle g_0g_1 \rangle$.
Consider the short exact sequence
\[
0 \longrightarrow R(-4,0)/(K:g_0g_1)  \stackrel{[g_0g_1]}{\longrightarrow}
R/K \longrightarrow R/J \longrightarrow 0.
\] 
{}From our choice of $\{g_0z,g_1w\}$ as a complete intersection,
it follows that the depth of $K$ is at least two. Since the maximal
minors of
\[
\left[
\begin{matrix}
g_0 & g_1  \\
-w & 0 \\
0 & -z
\end{matrix} \right]
\]
generate $K$, the second part of Theorem~\ref{hbst} implies that we
have a Hilbert-Burch resolution:
\[
\begin{aligned}
0 \to R(-2,-2)^2 &\xrightarrow{\left[
\begin{matrix}
g_0 & \!\!g_1  \\
-w & \!\!0 \\
0 & \!\!\!-z
\end{matrix} \right]} \!\!\!\!
\begin{array}{c}
R(0,-2) \\ \oplus \\ R(-2,-1)^2 \end{array}\!\!\!\!
\xrightarrow{\left[ \begin{matrix}
zw & \!\!g_0z & \!\!g_1w
\end{matrix}\right]} R \to R/K \to 0.
\end{aligned}
\]

Now we need to compute the free resolution of $(K : g_0g_1)$.  Observe
that both $z$ and $w$ are in the quotient and generate the entire
maximal ideal in $k[z,w]$.  Any element of $K$ has degree at least one
in these variables, and $g_0g_1$ has degree zero in these
variables. Thus, any element of $(K : g_0g_1)$ has degree at least one
in $z,w$, and so $z$ and $w$ must generate the whole ideal
quotient. In other words, $(K : g_0g_1) = \langle z,w \rangle$

Since $\{z,w\}$ is a regular sequence, the resolution of $R(-4,0)/(K :
g_0g_1)$ is a Koszul complex.  When we combine this with the above
Hilbert-Burch resolution for $R/K$ as in \eqref{beginnings}, we obtain
the diagram:

\vspace*{1.15in}

\begin{equation}
\label{rjcone}
\end{equation}

\vspace*{-1.85in}

\[
\quad\ \ \ \xymatrix@C=7pt{
 &&&& 0 \ar[d] \\
 0 \ar[r] & R(-4,-2) \ar[r] \ar[d] &  R(-4,-1)^2 \ar[r] \ar[d]
 & R(-4,-0) \ar[r] \ar[d] &
 R(-4,-0)/\langle z,w\rangle \ar[d]^{[g_0g_1]} \ar[r] & 0\\
 0 \ar[r] \ar@{--}[ur] & R(-2,-2)^2 \ar@{--}[ur] \ar[r] & 
{\begin{matrix} R(0,-2) \\ \oplus \\ R(-2,-1)^2 \end{matrix}}
\ar@{--}[ur] \ar[r] & 
  R \ar[r] \ar[dr] & R/K \ar[d]\ar[r] &0\\
 &&&&R/J \ar[d]\\
 &&&& 0 \\
}
\]

\noindent The mapping cone is constructed by taking the direct sum of
modules connected by dashed lines in the diagram.  This gives a free
resolution for $R/J$ as follows:
\[
0  \rightarrow  R(-4,-2) \rightarrow
\begin{array}{c}
R(-4,-1)^2\\
\oplus\\
R(-2,-2)^2 \\
\end{array}
 \rightarrow
\begin{array}{c}
 R(0,-2)\\
 \oplus \\
R(-4,0)\\
\oplus\\
R(-2,-1)^2 \\
\end{array}
\rightarrow R \rightarrow R/J \rightarrow 0.
\]
In fact, there is no overlap between the degree shifts,
and so in fact the mapping cone resolution is minimal.

We are now in a position to construct a mapping cone resolution for
$I$.  Rather surprisingly, it, too, turns out to be a minimal free
resolution.  Recall that for $R/I$, we have the exact sequence
\[
0 \longrightarrow R(-2,-1)/J \stackrel{[g_2l_2]
}{\longrightarrow} R/C \longrightarrow R/I \longrightarrow 0.
\]
where $C = \langle p_0,p_1\rangle = \langle g_0z,g_1w\rangle$ and $J =
(C:p_2) = (C:g_2l_2)$.  Thus, if we take the Koszul resolution of
$R/C$ and tensor the above resolution of $R/J$ with $R(-2,-1)$, we get
a mapping cone diagram similar to \eqref{rjcone} for $R/I$.  We
recommend that you write this out carefully.  You will see that there
is no overlap between degrees at consecutive steps of the resolution.
Hence, the mapping cone resolution is minimal.  Putting it all
together, we have proved the following theorem.

\begin{theorem}
\label{th:nongenres}
In the non-generic case, the minimal free resolution of the ideal $I =
 \langle p_0,p_1,p_2 \rangle$ is given by:
\[
0  \rightarrow  R(-6,-3) \rightarrow
\begin{array}{c}
R(-4,-3)^2\\
\oplus\\
R(-6,-2)^2 \\
\end{array}
 \rightarrow
\begin{array}{c}
 R(-6,-1)\\
 \oplus \\
R(-4,-2)^3\\
\oplus\\
R(-2,-3) \\
\end{array}
 \rightarrow
R(-2,-1) ^3
\rightarrow I \rightarrow 0.
\]
\end{theorem}

\subsection{The Generic Case}
Next, we tackle the generic case. It turns out that, just as in the
non-generic case, we can obtain the free resolution from an iterated
mapping cone construction.  We again make use of the explicit formulas
for the non-Koszul syzygies given in Section~6.2.  This will give
explicit generators for the ideal quotient.  By unwinding the internal
structure of these generators, we can use the same techniques to
construct the minimal free resolution.

\begin{theorem}
\label{th:genres}
In the generic case, the minimal free resolution of the ideal $I =
 \langle p_0,p_1,p_2 \rangle$ is given by:
\[
0  \rightarrow  R(-6,-3) \rightarrow
\begin{array}{c}
R(-4,-3)^3\\
\oplus\\
R(-6,-2)^2 \\
\end{array}
 \rightarrow
\begin{array}{c}
 R(-6,-1)\\
 \oplus \\
R(-4,-2)^3\\
\oplus\\
R(-3,-3)^2 \\
\end{array}
 \rightarrow
R(-2,-1) ^3
\rightarrow I \rightarrow 0.
\]
\end{theorem}

\begin{proof}
As in the non-generic case, we begin by choosing generators so that $I
= \langle p_0,p_1,p_2 \rangle$ with $C = \langle p_0,p_1 \rangle$ a
complete intersection. This is possible since two polynomials define a
complete intersection, unless they share a common factor.  On
$\mathbb{P}^1 \times \mathbb{P}^1$, such a common factor defines a
curve of degree $(a,b)$, which would meet $\mathbf{V}(p_2)$ in $2b+a$
points (see \cite[Example V.1.4.3]{h}), a contradiction to the
assumption that $I$ has no common zeros.  So there exists an exact
sequence
\begin{equation}
\label{rci}
0 \longrightarrow R(-2,-1)/(C : p_2) \stackrel{[p_2]}{\longrightarrow}
R/C \longrightarrow R/I \longrightarrow 0.
\end{equation}
Proposition~\ref{gensyz} gives the degrees of the generators of the
syzygy module $\mathrm{Syz}(p)$, and as in the non-generic case, this
gives the degrees of generators of $(C:p_2)$.  Hence $(C:p_2)$ has two
generators $p_0,p_1$ of degree $(2,1)$ corresponding to Koszul
syzygies, two generators $k_1,k_2$ of degree $(1,2)$, and a generator
$g$ of degree $(4,0)$.

Proposition~\ref{prop:deg33} gives explicit descriptions for the
generators $k_1,k_2$.  If we write
\begin{align*}
p_0 &= C_0x^2+D_0xy+E_0y^2\\
p_1 &= C_1x^2+D_1xy+E_1y^2,
\end{align*}
then we have
\[
k_1 = \det \left[\begin{matrix}
          C_0x+D_0y & E_0 \\ C_1x+D_1y & E_1
            \end{matrix}\right], \quad
k_2 = \det \left[\begin{matrix}
          C_0 & D_0x+ E_0y \\ C_1 & D_1x+ E_1y
            \end{matrix}\right].
\]

To understand $(C:p_2) = \langle p_0,p_1, k_1,k_2,g \rangle$, we will
use $K = \langle p_0,p_1, k_1,k_2 \rangle$ and $(K:g)$.  A check shows
that the $3\times3$ minors of the matrix:
\[
\phi = \left[\begin{matrix}
E_1 & -D_1 & C_1\\
E_0 & -D_0 & C_0\\
x & y & 0 \\
0 & x & y 
\end{matrix}\right]
\]
generate $K$.  Since the ideal generated by $\{p_0,p_1\}$ is a complete
intersection, $\{p_0,p_1\}$ is a regular sequence, hence the depth of
$K$ is at least two.  This means that we have a Hilbert-Burch
resolution, so $\phi$ is injective and we have a resolution
\[
\begin{aligned}
0 \rightarrow 
R(-2,-2)^3  \xrightarrow{ 
\phi } 
\begin{array}{c}
R(-1,-2)^2\\
\oplus\\
R(-2,-1)^2 \\
\end{array}
\xrightarrow{\left[\begin{matrix}
p_0 &\! p_1&\!\! -k_1 & \!\!-k_2
\end{matrix}\right]} R \rightarrow R/K \rightarrow 0
\end{aligned}
\]
Since $(C:p_2) = K+\langle g \rangle$, we get a mapping cone
resolution of $R/(C:p_2)$ (again, possibly non-minimal) from
resolutions of $R/K$ (given above) and $R(-4,0)/(K:g)$ (still to be
determined).  

To study $(K:g)$, recall from Proposition~\ref{prop:deg61}
that if we write 
\[
p_0 = A_0z+B_0w \quad\text{and}\quad p_1 = A_1z+B_1w,
\]
then $g= A_0B_1-A_1B_0$.  This makes it easy to see that $z,w \in
(K:g)$.  {}From here on, everything works as in the non-generic case:\
check that $(K:g) = \langle z,w \rangle$, and construct a mapping cone
resolution for $R/(C:p_2)$, which is forced to be minimal from degree
considerations.  By \eqref{rci}, the free resolutions of
$R(-2,-1)/(C:p_2)$ and $R/C$ give a mapping cone resolution for $R/I$,
which is also forced by degree considerations to be a minimal free
resolution.
\end{proof}

\subsection{Minimal Resolutions} One point we want to emphasize again is
that we were lucky that the resolutions resulting from the
 mapping cone construction were minimal.  In practice, one
usually has to ``prune'' the resolution before obtaining a minimal
one.  The paper \cite{mmr} gives some nice examples of how this is
done.  Finally, the fact that we had only
three generators need not imply that the resolution will be simple;
as Bruns shows in \cite{bruns}, essentially any resolution arises
as the resolution of a three generated ideal.

\subsection{Liason} The strategy of understanding $I = \langle
p_0,p_1,p_2\rangle$ by applying the mapping cone construction to $C =
\langle p_0,p_1\rangle$ and $J = (C : p_2)$ goes back to the classic paper
\cite{ps} of Peskine and Szpiro.  In general, ideals $I$ and $J$ of
codimension $c$ in a polynomial ring $R$ are said to be \emph{directly
linked} if there is a regular sequence $f_1,\dots,f_c$ contained in
$R$ such that
\begin{equation}
\label{liason}
J = (\langle f_1,\dots,f_c\rangle : I) \quad\text{and}\quad 
I = (\langle f_1,\dots,f_c\rangle : J).
\end{equation}
When $R/I$ is Cohen-Macaulay, then one can use the mapping cone
construction and duality to transform a free resolution of $R/I$ into
a free resolution of $R/J$.  See, for example, \cite[Exercise
21.23]{e}.  An introduction to \emph{linkage} (also called
\emph{liason}) can be found in \cite{m}.

In our situation, Bermejo and Gimenez observed that the symmetry of
\eqref{liason} fails; in particular (writing $C=\langle
p_0,p_1\rangle$) we have $J = (C : I)$, but a computer check shows $I
\ne (C : J)$. However, in the generic case, we
proved that $(K :
J) = (K:g) = \langle z,w\rangle$, and Bermejo and Gimenez noted that 
$K:\langle z,w\rangle = J$.
So we recover the symmetry \eqref{liason}, though not quite as liason
since $K$ is not a complete intersection.
While the usual notion of liason (using a complete intersection) 
is a bit too restrictive for our 
situation, there are more general notions of liason (see \cite{an}
or \cite{hu}) which may be useful in understanding this symmetry.

\section{The implicitization problem in geometric modeling}

A central problem in geometric modeling is to find the implicit
equations for a curve or surface defined by a rational map. For
surfaces, the two most common situations are the images of
parameterizations $\P^1 \times \P^1 \longrightarrow \P^3$ or $\P^2
\longrightarrow \P^3$. Surfaces of the first type are called {\it
tensor product surfaces}, and surfaces of the latter type are called
{\it triangular surfaces}.

The implicitization problem involves some interesting commutative
algebra.  For example, implicitization can be recast as a problem in
elimination theory, which allows one to use standard tools such as
Gr\"obner bases or resultants.  The surprise is that more
sophisticated tools from commutative algebra are also being used, and
syzygies play a leading role!  A detailed survey of this area appears
in \cite{c01}; we sketch some of the highlights below.

\subsection{Moving Curves and Surfaces}
In \cite{sc}, Sederberg and Chen introduced the method of moving
curves and surfaces.  For a curve parametrization $(a,b,c)$, a {\it moving
line} that follows the parametrization is a element of the syzygy
module on the generators of the ideal $I = \langle a,b,c\rangle$.  The
syzygy module of $I$ is free of rank two, and a Hilbert function
computation shows that if $a,b,c$ are homogeneous of degree $n$
without common factors, then there is an $n$-dimensional vector space
of moving lines of degree $n-1$.  Write each moving line as
\[
A_i(s,t)x+B_i(s,t)y+C_i(s,t)z
\]
where the $x,y,z$ are placeholders, representing the fact that
$A_i\cdot a + B_i \cdot b +C_i \cdot c = 0$. By collecting coefficients, 
we can write 
\[
A_i(s,t)x+B_i(s,t)y+C_i(s,t)z = \sum_{j=0}^{n-1}
L_{ij}(x,y,z)s^jt^{n-1-j}.
\]
A main theorem of \cite{csc} is that the determinant of the $n\times
n$ matrix of the $L_{ij}$ is a power of the implicit equation for the
image.

For a surface parametrization given by $(a,b,c,d)$, a moving plane
that follows the parametrization is an element of the syzygy module on
the generators of the ideal $I = \langle a,b,c,d\rangle$, and a moving
quadric that follows the parametrization is an element of the syzygy
module on the generators of $I^2$.  The moving surface method of
\cite{cgz} requires knowing that a syzygy of the form
\[
(c_1 a + c_2 b + c_3 c + c_4 d) a + (c_5 b + c_6 c + c_7 d) b + (c_8 c
+ c_9 d) c = 0
\]
comes from the Koszul complex when $a,b,c,d$ have no common zeros.  In
the case of $\P^2$ (i.e., when $a,b,c,d$ are homogeneous polynomials),
this is proved by observing that $a,b,c$ form a regular sequence, so
that every syzygy comes from the Koszul complex.  For $\P^1 \times
\P^1$ (i.e., when $a,b,c,d$ are bihomogeneous polynomials), the Koszul
complex is not exact in all bidegrees, but by vanishing theorems for
cohomology and arguments similar to those explained in our first sections, 
it can be seen that the sequence is exact in the bidegree of interest.

The main result of \cite{cgz} involves parameterizations without base
points.  When base points are allowed, \cite{bcd} shows that for
$\P^2$, the moving surface method of \cite{cgz} applies when the base
points are local complete intersections.  This is also true for
$\P^1\times\P^1$, by \cite{ahw}.  The proofs in \cite{bcd} use results
about the regularity of $I$ and $I^2$; in a similar way, the proofs in
\cite{ahw} use results about bigraded regularity.

Other directions of recent research are \cite{ccl}, which uses a
special case of the Serre conjecture to conclude that syzygy modules
are always free for affine surface parameterizations, and the use of
approximation complexes (\cite{bj}, \cite{bc}) to study
implicitization and moving surfaces.  The results above suggest that
the interaction between commutative algebra and implicitization is
unusually rich.

\subsection{Group Actions}
While the tools of commutative algebra can be more complicated in the
bigraded case, there are sometimes other tools that can be helpful.
In the case of a \emph{generic} tensor product surface of degree
$(2,1)$ in $\P^3$, the recent paper \cite{egl} uses group actions to
analyze the geometry of the resulting surface.  Here is a quick sketch
of their results.  Given four generic polynomials $p_0,p_1,p_2,p_3$ of
degree $(2,1)$, one can use group actions by $\mathrm{GL}(4)$ and
$\mathrm{GL}(2) \times \mathrm{GL}(2)$ to reduce to the case 
\[
p_0 =
x^2z,\ p_1 = (x-y)^2(z-w),\ p_2 = (x-Ay)^2(z-Bw),\ p_3 = y^2w,
\] 
where $A$
and $B$ are parameters determined by certain cross ratios.  Then 
Elkadi, Galligo and L\^{e}  show how to obtain the implicit
equation of the surface in terms of a certain $4\times4$ matrix, and
prove that the surface is singular along a twisted cubic.

In our situation, we have three polynomials $p_0,p_1,p_2$ of degree
$(2,1)$.  Similar to the previous paragraph,  there is the action
of $\mathrm{GL}(3)$ that replaces $p_0,p_1,p_2$ with linear
combinations, and second, the usual action of $\mathrm{GL}(2) \times
\mathrm{GL}(2)$ on $\P^1 \times \P^1$.  In our non-generic case, we
used these actions in Lemma~\ref{lem:cigens} to put $p_0,p_1,p_2$ in a
normal form.  

In our generic case, we know that $W(p) = \mathrm{Span}(p_0,p_1,p_2)$
meets $Y \simeq \P^2 \times \P^1$ in finitely many points.  Since
$Y \subset \P^5$ has degree three, the nicest case is when
$W(p)$ meets $Y$ in three distinct points.  When this happens, we can
pick $p_0, p_1, p_2$ to be these three points, and we get the
following normal form:
\begin{align*} 
p_0 &= x (x-Ay)z\\
p_1 &= y (y-Bx)w\\
p_2 &= (x+y)(x-Cy)(z+w).
\end{align*}
As one can easily guess, the resultant in Section
\ref{subs:resultants} is simply
\[ 
\mathrm{Res}_{(2,1)} = C(1+B)(BC-1)(1+A)(A-C)(AB-1).
\]
It is easy to check that the components of the syzygy
$C^{(1)}$ (resp.$C^{(2)}$) from Section~\ref{explicit} have a common
linear factor $Bx-y$ (resp. $-x+Ay$) if and only if $A =
BC+B+C-1$.  Thus, we are in the non-generic case precisely when this
equality holds.  However, it may happen that $W(p) =
\mathrm{Span}(p_0,p_1,p_2)$ meets $Y$ in a double point and a simple
point.  Consider
\begin{align*} 
p_0 &= x (x-Ay)z\\
p_1 &= y (y-Bx)z + x (x-Ay) w\\
p_2 &= (x+y)(x-Cy)(z+w),
\end{align*}
where $(1-AB)(A+1)(A-C)C \ne 0$.  In this case, one can show that
$W(p) \cap Y$ consists of only two points, $p_0$ and $p_2$, where
$p_0$ has multiplicity two.  One can also check that the line through
$p_0$ and $p_1$ lies in the tangent space to $Y$ at $p_0$, which shows
that the intersection $W(p) \cap Y$ is not transversal at $p_0$.

It follows that there is not a single ``normal form'' that covers all
cases in which $p_0, p_1, p_2$ have no common roots on $\P^1 \times
\P^1$, even in our generic case when $W(p)\cap Y$ is finite.

\section*{Acknowledgements} 
Evidence for this work was provided by many computations done using
Macaulay 2.  Macaulay 2 is available at the URL
\[
{\tt http://www.math.uiuc.edu/Macaulay2/}
\]
We also thank Isabel Bermejo and Philippe Gimenez for
helpful comments on the mapping cone strategy and liason.
Alicia Dickenstein was partially supported by UBACYT X042 and ANPCYT
03-6568, Argentina.  Hal Schenck was partially supported by NSF
DMS03-11142 and NSA MDA904-0301-0006. We are grateful to MSRI,
where some of this work took place.

\vskip .2in
\noindent\address{David Cox\\
Department of Mathematics and Computer Science, Amherst College\\
Amherst, MA 01002-5000\\}
{\tt dac@cs.amherst.edu}
\vskip .2in

\noindent\address{Alicia Dickenstein\\
Departamento de Matem\'atica, F.C.E.\ y N., Universidad de Buenos
Aires\\
Cuidad Universitaria--Pabell\'on I, 1428 Buenos Aires, Argentina\\}
{\tt alidick@dm.uba.ar}

\vskip .2in
\noindent\address{Hal Schenck\\
Department of Mathematics, Texas A{\&}M University\\
College Station, TX 77843\\}
{\tt schenck@math.tamu.edu}

\end{document}